\documentclass{article}
\usepackage{amsmath,amssymb}
\usepackage{cite}

\begin{document}

\begin{center}
{\bf New determinantal representations of the W-weighted Drazin inverse  over the
quaternion skew field}.\end{center}
\begin{center}{\bf Ivan Kyrchei}\end{center}\begin{center} Pidstrygach Institute for Applied Problems of Mechanics and Mathematics of NAS of Ukraine,
 Ukraine, kyrchei@online.ua\end{center}

\begin{abstract}
Within the framework of the theory of the column and row determinants, we obtain new determinantal representations of the W-weighted Drazin inverse over the quaternion skew field. We give determinantal representations of the W-weighted Drazin inverse by using previously introduced  determinantal representations of the Drazin inverse, the Moore-Penrose inverse, and the limit representations of the W-weighted Drazin inverse in some special case.
\end{abstract}

\section{Introduction}
\newtheorem{Corollary}{Corollary}[section]
\newtheorem{theorem}{Theorem}[section]
\newtheorem{lemma}{Lemma}[section]
\newtheorem{definition}{Definition}[section]
\newtheorem{remark}{Remark}[section]
\newcommand{\rank}{\mathop{\rm rank}\nolimits}
\newtheorem{proposition}{Proposition}[section]

Let ${\rm
{\mathbb{R}}}$ and ${\rm
{\mathbb{C}}}$ be the real and complex number fields, respectively.
Throughout the paper, we denote  the set of all $m\times n$ matrices over the
quaternion algebra
\[{\rm {\mathbb{H}}}=\{a_{0}+a_{1}i+a_{2}j+a_{3}k\,
|\,i^{2}=j^{2}=k^{2}=-1,\, a_{0}, a_{1}, a_{2}, a_{3}\in{\rm
{\mathbb{R}}}\}\]
by ${\rm {\mathbb{H}}}^{m\times n}$, and by ${\rm {\mathbb{H}}}^{m\times n}_{r}$ the set of all $m\times n$ matrices over $\mathbb{H}$ with a rank $r$. Let ${\rm M}\left( {n,{\rm {\mathbb{H}}}} \right)$ be the
ring of $n\times n$ quaternion matrices and $ {\bf I}$ be the identity matrix with the appropriate size. For ${\rm {\bf A}}
 \in {\rm {\mathbb{H}}}^{n\times m}$, we denote by ${\rm {\bf A}}^{ *}$, $\rank {\bf A}$ the conjugate transpose (Hermitian adjoint) matrix and the rank
of ${\rm {\bf A}}$.
 The matrix ${\rm {\bf A}} = \left( {a_{ij}}  \right) \in {\rm
{\mathbb{H}}}^{n\times n}$ is Hermitian if ${\rm {\bf
A}}^{ *}  = {\rm {\bf A}}$.

The  definitions of the generalized inverse matrices may be  extended  to quaternion matrices.

The Moore-Penrose inverse of ${\bf A}\in{\rm {\mathbb{H}}}^{m\times n}$, denoted by ${\bf A}^{\dagger}$, is the unique matrix ${\bf X}\in{\rm {\mathbb{H}}}^{n\times m}$ satisfying the following equations,
 \begin{gather}\label{eq1:MP}  {\rm {\bf A}}{\bf X}
{\rm {\bf A}} = {\rm {\bf A}}; \\
                                  {\bf X} {\rm {\bf
A}}{\bf X}  = {\bf X};\\
                                  \left( {\rm {\bf A}}{\bf X} \right)^{ *}  = {\rm
{\bf A}}{\bf X}; \\
                                \left( {{\bf X} {\rm {\bf A}}} \right)^{ *}  ={\bf X} {\rm {\bf A}}. \end{gather}
 For ${\bf A}\in{\rm {\mathbb{H}}}^{n\times n}$ with $k = Ind\,{\bf A}$ the smallest positive number such that $\rank {\bf A}^{k+1}=\rank {\bf A}^{k}$ the
Drazin inverse of ${\bf A}$, denoted by ${\bf A}^{D}$, is defined to be the unique matrix ${\bf X}$ that satisfying (1.2)
and the following equations,
                                \begin{gather}\label{eq:Dr}
                                   {\rm {\bf A}}{\bf X}  = {\rm
{\bf A}}{\bf X}; \\
                                   {\rm {\bf A}}^{k+1}{\rm
{\bf X}}={\rm {\bf
         A}}^{k}. \end{gather}
 In particular, when $Ind{\kern 1pt} {\rm {\bf A}}=1$,
then the matrix ${\rm {\bf X}}$  is called
the group inverse and is denoted by ${\rm {\bf X}}={\rm {\bf
A}}^{g }$.

If $Ind{\kern 1pt} {\rm {\bf A}}=0$, then ${\rm {\bf A}}$ is
nonsingular, and ${\rm {\bf A}}^{D}\equiv {\bf A}^{\dagger}= {\rm {\bf A}}^{-1}$.

Cline and Greville \cite{cl} extended the Drazin inverse of square matrix to rectangular matrix, which can be generalized to the quaternion algebra as follows. For ${\bf A}\in{\rm {\mathbb{H}}}^{m\times n}$ and ${\bf W}\in{\rm {\mathbb{H}}}^{n\times m}$, the W-weighted Drazin inverse of ${\bf A}$ with respect to ${\bf W}$, denoted by ${\bf A}_{d,W}$, is the unique solution to equations,
\begin{gather}\label{eq:WDr}
  ({{\bf A}{\bf W}})^{k+1}{\rm
{\bf X}}{\bf W}=({\rm {\bf
         A}}{\bf W})^{k};\\
  {\rm {\bf X}}{\bf W}{\rm {\bf A}}{\bf W}{\rm {\bf X}}={\rm {\bf X}};\\
   {\rm
{\bf A}}{\bf W}{\rm {\bf X}}={\rm {\bf X}}{\bf W}{\rm {\bf A}},
\end{gather}
where $k= {\rm max}\{Ind({\bf A}{\bf W}), Ind({\bf W}{\bf A})\}$.
It is denoted by ${\rm {\bf X}}={\rm {\bf A}}_{d,{\bf W}}$.
The properties of the complex W-weighted Drazin inverse can be found in \cite{cl,wei1,wei2,wei4,ze,mos}.
These properties can be generalized to ${\mathbb{H}}$.
If ${\bf A}\in {\mathbb{H}}^{m\times n}$ with respect to ${\bf W}\in {\mathbb{H}}^{n\times m}$ and $k= {\rm max}\{Ind({\bf A}{\bf W}), Ind({\bf W}{\bf A})\}$, then
   \begin{equation}\label{eq:WDr_dr}
{\rm {\bf A}}_{d,{\bf W}}={\bf A}\left(({\bf W}{\bf A})^{D}) \right)^{2}=\left(({\bf A}{\bf W})^{D}) \right)^{2}{\bf A},
\end{equation}
   \begin{equation}
\label{eq:ADW}
{\rm {\bf A}}_{d,{\bf W}}{\bf W}=({\bf W}{\bf A})^{D}),  {\bf W}{\rm {\bf A}}_{d,{\bf W}}=({\bf A}{\bf W})^{D}.
\end{equation}
The problem of determinantal representation of generalized inverse matrices only recently begun to be decided through the theory of the column-row determinants introduced in \cite{ky1,ky2}.
The theory of row and column determinants  develops
the classical approach to a definition of a determinant, as alternating sum
of products of the entries of matrix but with a predetermined order of factors in each
terms of the determinant. A determinant of a quadratic matrix with noncommutative elements is often called the noncommutative determinant. Unlike other known noncommutative determinants such as determinants of Dieudonn\'{e} \cite{di}, Study \cite{stu}, Moore \cite{mo,dy}, Chen  \cite{ch}, quasideterminants of Gelfand-Retakh  \cite{ge}, the double determinant  built on the theory of the column-row determinants has properties similar to a usual determinant, in particular it can be expand  along arbitrary rows and columns. This property is necessary for determinantal representations of an inverse and generalized inverse matrices.
Determinantal representations of the Moore-Penrose inverse and the Drazin inverse  over the quaternion skew-field have been obtained in \cite{ ky3,ky33} and \cite{ky4}, respectively. Determinantal representations of an outer inverse ${\bf A}_ {T, S}^{(2)}$ is introduced in \cite{song1,song3} using the column-row determinants as well. Recall that an outer inverse of a matrix ${\bf A}$ over complex field with prescribed range
space $T$ and null space $S$ is a solution of (1.2) with restrictions,
\[\mathcal{R}({\bf X})=T,\,\,\,\mathcal{N}({\bf X})=S.\]
   Within the framework of the theory of the column-row determinants Song \cite{song2} also gave a determinantal representation W-weighted Drazin inverse over the quaternion skew-field using a characterization of the W-weighted Drazin
inverse by an outer inverse ${\bf A}_ {T, S}^{(2)}$. But  in obtaining of this determinantal representation    is used auxiliary matrices which different from ${\bf A}$  or its powers.
In this paper we  obtain    determinantal representations of the W-weighted Drazin inverse of ${\bf A}\in{\rm {\mathbb{H}}}^{m\times n}$ with respect to  ${\bf W}\in{\rm {\mathbb{H}}}^{n\times m}$ by using only their entries.

The paper is organized as follows. We start with some
basic concepts and results from the theory of  the row and column determinants and give the determinantal representations of the inverse, the Moore-Penrose inverse, and the Drazin inverse over the quaternion skew field  in Section 2.  In Section 3, we obtain determinantal representations of the W-weighted Drazin inverse by using introduced above determinantal representations of the Drazin inverse, the Moore-Penrose inverse, and the limit representations of the W-weighted Drazin inverse in some special case. In Section 4, we show
a numerical example to illustrate the main result.
\section{Elements of the theory of the column and row determinants}
 For a quadratic matrix ${\rm {\bf A}}=(a_{ij}) \in {\rm
M}\left( {n,{\mathbb{H}}} \right)$ we define $n$ row determinants and $n$ column determinants as follows.
Suppose $S_{n}$ is the symmetric group on the set $I_{n}=\{1,\ldots,n\}$.
\begin{definition}
 The $i$th row determinant of ${\rm {\bf A}}=(a_{ij}) \in {\rm
M}\left( {n,{\mathbb{H}}} \right)$ is defined  for all $i = \overline{1,n} $
by putting
 \begin{gather*}{\rm{rdet}}_{ i} {\rm {\bf A}} =
{\sum\limits_{\sigma \in S_{n}} {\left( { - 1} \right)^{n - r}{a_{i{\kern
1pt} i_{k_{1}}} } {a_{i_{k_{1}}   i_{k_{1} + 1}}} \ldots } } {a_{i_{k_{1}
+ l_{1}}
 i}}  \ldots  {a_{i_{k_{r}}  i_{k_{r} + 1}}}
\ldots  {a_{i_{k_{r} + l_{r}}  i_{k_{r}} }},\\
\sigma = \left(
{i\,i_{k_{1}}  i_{k_{1} + 1} \ldots i_{k_{1} + l_{1}} } \right)\left(
{i_{k_{2}}  i_{k_{2} + 1} \ldots i_{k_{2} + l_{2}} } \right)\ldots \left(
{i_{k_{r}}  i_{k_{r} + 1} \ldots i_{k_{r} + l_{r}} } \right),\end{gather*}
with
conditions $i_{k_{2}} < i_{k_{3}}  < \ldots < i_{k_{r}}$ and $i_{k_{t}}  <
i_{k_{t} + s} $ for $t = \overline{2,r} $ and $s =\overline{1,l_{t}} $.
\end{definition}
\begin{definition}
The $j$th column determinant
 of ${\rm {\bf
A}}=(a_{ij}) \in {\rm M}\left( {n,{\mathbb{H}}} \right)$ is defined for
all $j =\overline{1,n} $ by putting
 \begin{gather*}{\rm{cdet}} _{{j}}\, {\rm {\bf A}} =
{{\sum\limits_{\tau \in S_{n}} {\left( { - 1} \right)^{n - r}a_{j_{k_{r}}
j_{k_{r} + l_{r}} } \ldots a_{j_{k_{r} + 1} i_{k_{r}} }  \ldots } }a_{j\,
j_{k_{1} + l_{1}} }  \ldots  a_{ j_{k_{1} + 1} j_{k_{1}} }a_{j_{k_{1}}
j}},\\
\tau =
\left( {j_{k_{r} + l_{r}}  \ldots j_{k_{r} + 1} j_{k_{r}} } \right)\ldots
\left( {j_{k_{2} + l_{2}}  \ldots j_{k_{2} + 1} j_{k_{2}} } \right){\kern
1pt} \left( {j_{k_{1} + l_{1}}  \ldots j_{k_{1} + 1} j_{k_{1} } j}
\right), \end{gather*}
\noindent with conditions, $j_{k_{2}}  < j_{k_{3}}  < \ldots <
j_{k_{r}} $ and $j_{k_{t}}  < j_{k_{t} + s} $ for  $t = \overline{2,r} $
and $s = \overline{1,l_{t}}  $.
\end{definition}
Suppose ${\rm {\bf A}}_{}^{i{\kern 1pt} j} $ denotes the submatrix of
${\rm {\bf A}}$ obtained by deleting both the $i$th row and the $j$th
column. Let ${\rm {\bf a}}_{.j} $ be the $j$th column and ${\rm {\bf
a}}_{i.} $ be the $i$th row of ${\rm {\bf A}}$. Suppose ${\rm {\bf
A}}_{.j} \left( {{\rm {\bf b}}} \right)$ denotes the matrix obtained from
${\rm {\bf A}}$ by replacing its $j$th column with the column ${\rm {\bf
b}}$, and ${\rm {\bf A}}_{i.} \left( {{\rm {\bf b}}} \right)$ denotes the
matrix obtained from ${\rm {\bf A}}$ by replacing its $i$th row with the
row ${\rm {\bf b}}$.

The following theorem has a key value in the theory of the column and row
determinants.
\begin{theorem} \cite{ky1}\label{theorem:
determinant of hermitian matrix} If ${\rm {\bf A}} = \left( {a_{ij}}
\right) \in {\rm M}\left( {n,{\rm {\mathbb{H}}}} \right)$ is Hermitian,
then ${\rm{rdet}} _{1} {\rm {\bf A}} = \cdots = {\rm{rdet}} _{n} {\rm {\bf
A}} = {\rm{cdet}} _{1} {\rm {\bf A}} = \cdots = {\rm{cdet}} _{n} {\rm {\bf
A}} \in {\rm {\mathbb{R}}}.$
\end{theorem}
 Since all  column and  row determinants of a
Hermitian matrix over ${\rm {\mathbb{H}}}$ are equal, we can define the
determinant of a  Hermitian matrix ${\rm {\bf A}}\in {\rm M}\left( {n,{\rm
{\mathbb{H}}}} \right)$. By definition, we put
$\det {\rm {\bf A}}: = {\rm{rdet}}_{{i}}\,
{\rm {\bf A}} = {\rm{cdet}} _{{i}}\, {\rm {\bf A}}, $
 for all $i =\overline{1,n}$.

The determinant of a Hermitian matrix has properties similar to a usual determinant. They are completely explored
in
 \cite{ky1, ky2}
 by its row and
column determinants. They can be summarized by the following
theorems.
\begin{theorem}\label{theorem:row_combin} If the $i$th row of
a Hermitian matrix ${\rm {\bf A}}\in {\rm M}\left( {n,{\rm
{\mathbb{H}}}} \right)$ is replaced with a left linear combination
of its other rows, i.e. ${\rm {\bf a}}_{i.} = c_{1} {\rm {\bf
a}}_{i_{1} .} + \ldots + c_{k}  {\rm {\bf a}}_{i_{k} .}$, where $
c_{l} \in {{\rm {\mathbb{H}}}}$ for all $ l = \overline{1, k}$ and
$\{i,i_{l}\}\subset I_{n} $, then
\[
 {\rm{rdet}}_{i}\, {\rm {\bf A}}_{i \, .} \left(
{c_{1} {\rm {\bf a}}_{i_{1} .} + \ldots + c_{k} {\rm {\bf
a}}_{i_{k} .}}  \right) = {\rm{cdet}} _{i}\, {\rm {\bf A}}_{i\, .}
\left( {c_{1}
 {\rm {\bf a}}_{i_{1} .} + \ldots + c_{k} {\rm {\bf
a}}_{i_{k} .}}  \right) = 0.
\]
\end{theorem}
\begin{theorem}\label{theorem:colum_combin} If the $j$th column of
 a Hermitian matrix ${\rm {\bf A}}\in
{\rm M}\left( {n,{\rm {\mathbb{H}}}} \right)$   is replaced with a
right linear combination of its other columns, i.e. ${\rm {\bf
a}}_{.j} = {\rm {\bf a}}_{.j_{1}}   c_{1} + \ldots + {\rm {\bf
a}}_{.j_{k}} c_{k} $, where $c_{l} \in{{\rm {\mathbb{H}}}}$ for
all $ l = \overline{1, k}$ and $\{j,j_{l}\}\subset J_{n}$, then
 \[{\rm{cdet}} _{j}\, {\rm {\bf A}}_{.j}
\left( {{\rm {\bf a}}_{.j_{1}} c_{1} + \ldots + {\rm {\bf
a}}_{.j_{k}}c_{k}} \right) ={\rm{rdet}} _{j} \,{\rm {\bf A}}_{.j}
\left( {{\rm {\bf a}}_{.j_{1}}  c_{1} + \ldots + {\rm {\bf
a}}_{.j_{k}}  c_{k}} \right) = 0.
\]
\end{theorem}
  The determinant of a Hermitian matrix also has a property of expansion along arbitrary rows and columns using  row and column determinants of submatrices. So, we were able to get determinantal representations of an inverse and generalized inverse matrices as follows.
\begin{theorem}\cite{ky1} \label{inver_her} If for a Hermitian matrix ${\rm {\bf A}}\in {\rm M}\left( {n,{\rm
{\mathbb{H}}}} \right)$,
\[\det {\rm {\bf A}} \ne 0,\]
then there exist a unique right inverse  matrix $(R{\rm {\bf A}})^{ - 1}$
and a unique left inverse matrix $(L{\rm {\bf A}})^{ - 1}$ of a
nonsingular
 ${\rm {\bf A}}$, where $\left( {R{\rm {\bf A}}} \right)^{ - 1} = \left( {L{\rm {\bf A}}}
\right)^{ - 1} = :{\rm {\bf A}}^{ - 1}$, and the right  and left inverse
matrices possess the following determinantal representations
\begin{equation}\label{eq:inver_her_R}
  \left( {R{\rm {\bf A}}} \right)^{ - 1} = {\frac{{1}}{{\det {\rm
{\bf A}}}}}
\begin{pmatrix}
  R_{11} & R_{21} & \cdots & R_{n1}\\
  R_{12} & R_{22} & \cdots & R_{n2}\\
  \cdots & \cdots & \cdots& \cdots\\
  R_{1n} & R_{2n} & \cdots & R_{nn}
\end{pmatrix},
\end{equation}
\begin{equation}\label{eq:inver_her_L}
  \left( {L{\rm {\bf A}}} \right)^{ - 1} = {\frac{{1}}{{\det {\rm
{\bf A}}}}}
\begin{pmatrix}
  L_{11} & L_{21} & \cdots & L_{n1} \\
  L_{12} & L_{22} & \cdots & L_{n2} \\
  \cdots & \cdots & \cdots & \cdots \\
  L_{1n} & L_{2n} & \cdots & L_{nn}
\end{pmatrix},
\end{equation}
where $ \det {\rm {\bf A}} ={\sum\limits_{j = 1}^{n} {{a_{i{\kern
1pt} j} \cdot R_{i{\kern 1pt} j} } }}= {{\sum\limits_{i = 1}^{n}
{L_{i{\kern 1pt} j} \cdot a_{i{\kern 1pt} j}} }}$,
\begin{gather*}
 R_{i{\kern 1pt} j} = {\left\{ {{\begin{array}{*{20}c}
  - {\rm{rdet}}_{{j}}\, {\rm {\bf A}}_{{.{\kern 1pt} j}}^{{i{\kern 1pt} i}} \left( {{\rm
{\bf a}}_{{.{\kern 1pt} {\kern 1pt} i}}}  \right),& {i \ne j},
\hfill \\
 {\rm{rdet}} _{{k}}\, {\rm {\bf A}}^{{i{\kern 1pt} i}},&{i = j},
\hfill \\
\end{array}} } \right.}\,\,\,\,
L_{i{\kern 1pt} j} = {\left\{ {\begin{array}{*{20}c}
 -{\rm{cdet}} _{i}\, {\rm {\bf A}}_{i{\kern 1pt} .}^{j{\kern 1pt}j} \left( {{\rm {\bf a}}_{j{\kern 1pt}. } }\right),& {i \ne
j},\\
 {\rm{cdet}} _{k}\, {\rm {\bf A}}^{j\, j},& {i = j},
\\
\end{array} }\right.}
\end{gather*}
 and ${\rm {\bf A}}_{.{\kern 1pt} j}^{i{\kern 1pt} i} \left(
{{\rm {\bf a}}_{.{\kern 1pt} {\kern 1pt} i}}  \right)$ is obtained from
${\rm {\bf A}}$   by both replacing the $j$th column with the $i$th column and deleting  the $i$th row and column, ${\rm {\bf A}}_{i{\kern 1pt} .}^{jj} \left( {{\rm {\bf
a}}_{j{\kern 1pt} .} } \right)$ is obtained
 by both replacing the $i$th row with the $j$th row and
deleting  the $j$th row and  column, respectively, $I_{n} =  {\left\{ {1,\ldots ,n} \right\}}$, $k = \min {\left\{ {I_{n}}
\right.} \setminus {\left. {\{i\}} \right\}}$
 for all $ i,j =
\overline{1,n}$.
\end{theorem}
We shall use the following notations. Let $\alpha : = \left\{
{\alpha _{1} ,\ldots ,\alpha _{k}} \right\} \subseteq {\left\{
{1,\ldots ,m} \right\}}$ and $\beta : = \left\{ {\beta _{1}
,\ldots ,\beta _{k}} \right\} \subseteq {\left\{ {1,\ldots ,n}
\right\}}$ be subsets of the order $1 \le k \le \min {\left\{
{m,n} \right\}}$. By ${\rm {\bf A}}_{\beta} ^{\alpha} $ denote the
submatrix of ${\rm {\bf A}}$ determined by the rows indexed by
$\alpha$ and the columns indexed by $\beta$. Then ${\rm {\bf
A}}{\kern 1pt}_{\alpha} ^{\alpha}$ denotes the principal submatrix
determined by the rows and columns indexed by $\alpha$.
 If ${\rm {\bf A}} \in {\rm
M}\left( {n,{\rm {\mathbb{H}}}} \right)$ is Hermitian, then by
${\left| {{\rm {\bf A}}_{\alpha} ^{\alpha} } \right|}$ denote the
corresponding principal minor of $\det {\rm {\bf A}}$.
 For $1 \leq k\leq n$, the collection of strictly
increasing sequences of $k$ integers chosen from $\left\{
{1,\ldots ,n} \right\}$ is denoted by $\textsl{L}_{ k,
n}: = {\left\{ {\,\alpha :\alpha = \left( {\alpha _{1} ,\ldots
,\alpha _{k}} \right),\,{\kern 1pt} 1 \le \alpha _{1} \le \ldots
\le \alpha _{k} \le n} \right\}}$.  For fixed $i \in \alpha $ and $j \in
\beta $, let $I_{r,\,m} {\left\{ {i} \right\}}: = {\left\{
{\,\alpha :\alpha \in L_{r,m} ,i \in \alpha}  \right\}}{\rm ,}
\quad J_{r,\,n} {\left\{ {j} \right\}}: = {\left\{ {\,\beta :\beta
\in L_{r,n} ,j \in \beta}  \right\}}$.

 Denote by ${\rm {\bf a}}_{.j}^{*} $ and ${\rm {\bf
a}}_{i.}^{*} $ the $j$th column  and the $i$th row of  ${\rm
{\bf A}}^{*} $ and by ${\rm {\bf a}}_{.j}^{(m)} $ and ${\rm {\bf
a}}_{i.}^{(m)} $ the $j$th column  and the $i$th row of  ${\rm
{\bf A}}^{m} $, respectively.

The following theorem give determinantal representations of the Moore-Penrose inverse over the quaternion skew field $\mathbb{H}$.
\begin{theorem} \cite{ky3}\label{theor:det_repr_MP}
If ${\rm {\bf A}} \in {\rm {\mathbb{H}}}_{r}^{m\times n} $, then
the Moore-Penrose inverse  ${\rm {\bf A}}^{ +} = \left( {a_{ij}^{
+} } \right) \in {\rm {\mathbb{H}}}_{}^{n\times m} $ possess the
following determinantal representations:
\begin{equation}
\label{eq:det_repr_A*A}
 a_{ij}^{ +}  = {\frac{{{\sum\limits_{\beta
\in J_{r,\,n} {\left\{ {i} \right\}}} {{\rm{cdet}} _{i} \left(
{\left( {{\rm {\bf A}}^{ *} {\rm {\bf A}}} \right)_{\,. \,i}
\left( {{\rm {\bf a}}_{.j}^{ *} }  \right)} \right){\kern 1pt}
{\kern 1pt} _{\beta} ^{\beta} } } }}{{{\sum\limits_{\beta \in
J_{r,\,\,n}} {{\left| {\left( {{\rm {\bf A}}^{ *} {\rm {\bf A}}}
\right){\kern 1pt} _{\beta} ^{\beta} }  \right|}}} }}},
\end{equation}
or
\begin{equation}
\label{eq:det_repr_AA*} a_{ij}^{ +}  =
{\frac{{{\sum\limits_{\alpha \in I_{r,m} {\left\{ {j} \right\}}}
{{\rm{rdet}} _{j} \left( {({\rm {\bf A}}{\rm {\bf A}}^{ *}
)_{j\,.\,} ({\rm {\bf a}}_{i.\,}^{ *} )} \right)\,_{\alpha}
^{\alpha} } }}}{{{\sum\limits_{\alpha \in I_{r,\,m}}  {{\left|
{\left( {{\rm {\bf A}}{\rm {\bf A}}^{ *} } \right){\kern 1pt}
_{\alpha} ^{\alpha} } \right|}}} }}}.
\end{equation}
for all $i = \overline{1, n} $, $j =\overline{1, m} $.
\end{theorem}
\begin{proposition}\label{theor:repr_Dr} \cite{ca1} If $Ind({\bf A}) = k$, then
${\bf A}^{D}={\bf A}^{k}({\bf A}^{2k+1})^{+}{\bf A}^{k}
$.
\end{proposition}
Using  the determinantal representations of the Moore-Penrose inverse (\ref{eq:det_repr_A*A}) and (\ref{eq:det_repr_AA*}), and  Proposition \ref{theor:repr_Dr} we have obtained the following determinantal representations of the Drazin inverse for an arbitrary square matrix over $\mathbb{H}$.
Denote  by $\hat{{\rm {\bf a}}}_{.s}$ and $\check{{{\rm {\bf a}}}}_{t.}$ the $s$th column of $
({\bf A}^{ 2k+1})^{*} {\bf A}^{k}=:\hat{{\rm {\bf A}}}=
(\hat{a}_{ij})\in {\mathbb{H}}^{n\times n}$
and  the $t$th row of $
{\bf A}^{k}({\bf A}^{ 2k+1})^{*} =:\check{{\rm {\bf A}}}=
(\check{a}_{ij})\in {\mathbb{H}}^{n\times n}$, respectively,  for all $s,t=\overline
{1,n}$.
\begin{theorem} \cite{ky4}\label{theor:det_rep_draz} If ${\rm {\bf A}} \in {\rm M}\left( {n, {\mathbb{H}}}\right)$  with
$ Ind{\kern 1pt} {\rm {\bf A}}=k$ and $\rank{\rm {\bf A}}^{k+1} =
\rank{\rm {\bf A}}^{k} = r$, then the Drazin inverse  ${\rm {\bf A}}^{ D} $ possess the
determinantal representations
\begin{equation}
\label{eq:cdet_draz} a_{ij} ^{D}=
 {\frac{{ \sum\limits_{t = 1}^{n} {a}_{it}^{(k)}   {\sum\limits_{\beta \in J_{r,\,n} {\left\{ {t}
\right\}}} {{\rm{cdet}} _{t} \left( {\left({\bf A}^{ 2k+1} \right)^{*}\left({\bf A}^{ 2k+1} \right)_{. \,t} \left( \hat{{\rm {\bf a}}}_{.\,j}
\right)} \right){\kern 1pt}  _{\beta} ^{\beta} } }
}}{{{\sum\limits_{\beta \in J_{r,\,n}} {{\left| {\left({\bf A}^{ 2k+1} \right)^{*}\left({\bf A}^{ 2k+1} \right){\kern 1pt} _{\beta} ^{\beta}
}  \right|}}} }}}
\end{equation}
and
\begin{equation}
\label{eq:rdet_draz} a_{ij} ^{D}=
{\frac{\sum\limits_{s = 1}^{n}\left({{\sum\limits_{\alpha \in I_{r,\,n} {\left\{ {s}
\right\}}} {{\rm{rdet}} _{s} \left( {\left( { {\bf A}^{ 2k+1} \left({\bf A}^{ 2k+1} \right)^{*}
} \right)_{\,. s} (\check{{\rm {\bf a}}}_{i\,.})} \right) {\kern 1pt} _{\alpha} ^{\alpha} } }
}\right){a}_{sj}^{(k)}}{{{\sum\limits_{\alpha \in I_{r,\,n}} {{\left| {\left( { {\bf A}^{ 2k+1} \left({\bf A}^{ 2k+1} \right)^{*}
} \right){\kern 1pt} _{\alpha} ^{\alpha}
}  \right|}}} }}}
\end{equation}
\end{theorem}
In the special case, when ${\rm {\bf A}} \in {\rm M}\left( {n, {\mathbb{H}}}\right)$ is Hermitian, we can obtain  simpler determinantal
representations of the Drazin inverse.
\begin{theorem} \cite{ky4}\label{theor:det_rep_draz_her}
If ${\rm {\bf A}} \in {\rm M}\left( {n, {\mathbb{H}}}\right)$ is Hermitian with
$ Ind{\kern 1pt} {\rm {\bf A}}=k$ and $\rank{\rm {\bf A}}^{k+1} =
\rank{\rm {\bf A}}^{k} = r$, then the Drazin inverse ${\rm {\bf
A}}^{D} = \left( {a_{ij}^{D} } \right) \in {\rm
{\mathbb{H}}}_{}^{n\times n} $ possess the following determinantal
representations:
\begin{equation}
\label{eq:dr_rep_cdet} a_{ij}^{D}  = {\frac{{{\sum\limits_{\beta
\in J_{r,\,n} {\left\{ {i} \right\}}} {{\rm{cdet}} _{i} \left(
{\left( {{\rm {\bf A}}^{k+1}} \right)_{\,. \,i} \left( {{\rm {\bf
a}}_{.j}^{ k} }  \right)} \right){\kern 1pt} {\kern 1pt} _{\beta}
^{\beta} } } }}{{{\sum\limits_{\beta \in J_{r,\,\,n}} {{\left|
{\left( {{\rm {\bf A}}^{k+1}} \right){\kern 1pt} _{\beta} ^{\beta}
}  \right|}}} }}},
\end{equation}
or
\begin{equation}
\label{eq:dr_rep_rdet} a_{ij}^{D}  = {\frac{{{\sum\limits_{\alpha
\in I_{r,n} {\left\{ {j} \right\}}} {{\rm{rdet}} _{j} \left(
{({\rm {\bf A}}^{ k+1} )_{j\,.\,} ({\rm {\bf a}}_{i.\,}^{ (k)} )}
\right)\,_{\alpha} ^{\alpha} } }}}{{{\sum\limits_{\alpha \in
I_{r,\,n}}  {{\left| {\left( {{\rm {\bf A}}^{k+1} } \right){\kern
1pt}  _{\alpha} ^{\alpha} } \right|}}} }}}.
\end{equation}
\end{theorem}
Note that \textbf{the determinantal rank} of  ${\bf A}\in  {\mathbb{C}}^{m\times n}$ can be obtained as the largest order of a  non-zero principle minor in the Hermitian matrices ${\bf A}^{*}{\bf A}$ or ${\bf A}{\bf A}^{*}$.

We shall also need the following facts about the eigenvalues of quaternion matrices.
Due to the noncommutativity of quaternions, there are two
types of eigenvalues.
A quaternion $\lambda$ is said to be a right eigenvalue
of ${\rm {\bf A}} \in {\rm M}\left( {n,{\rm {\mathbb{H}}}}
\right)$ if ${\rm {\bf A}} \cdot {\rm {\bf x}} = {\rm
{\bf x}} \cdot \lambda $ for some nonzero quaternion column-vector
${\rm {\bf x}}$ with quaternion components. Similarly $\lambda$ is a left eigenvalue if ${\rm
{\bf A}} \cdot {\rm {\bf x}} = \lambda \cdot {\rm {\bf x}}$ for some nonzero quaternion column-vector
${\rm {\bf x}}$ with quaternion components.
The theory on the left eigenvalues of quaternion matrices has been
investigated in particular in \cite{hu, so, wo}. The theory on the
right eigenvalues of quaternion matrices is more developed. In
particular we note  \cite{br,ma,ba,dra,zh,far}.
\begin{proposition}\cite{zh}
Let ${\rm {\bf A}} \in {\rm M}\left( {n,{\rm {\mathbb{H}}}}
\right)$ is Hermitian. Then ${\rm {\bf A}}$ has exactly $n$ real
right eigenvalues.
\end{proposition}
Right and left eigenvalues are in general unrelated \cite{fa}, but it is not  for Hermitian matrices.
Suppose ${\rm {\bf A}} \in {\rm M}\left( {n, {\mathbb{H}}}\right)$
is Hermitian and $\lambda \in {\rm {\mathbb {R}}}$ is its right
eigenvalue, then ${\rm {\bf A}} \cdot {\rm {\bf x}} = {\rm {\bf
x}} \cdot \lambda = \lambda \cdot {\rm {\bf x}}$. This means that
all right eigenvalues of a Hermitian matrix are its left
eigenvalues as well. For real left eigenvalues, $\lambda \in {\rm
{\mathbb {R}}}$, the matrix $\lambda {\rm {\bf I}} - {\rm {\bf
A}}$ is Hermitian.
\begin{definition}
If $t \in {\rm {\mathbb {R}}}$, then for a Hermitian matrix ${\rm {\bf A}} \in {\rm M}\left( {n,{\rm {\mathbb{H}}}}
\right)$ the polynomial $p_{{\rm {\bf A}}}\left( {t} \right) =
\det \left( {t{\rm {\bf I}} - {\rm {\bf A}}} \right)$ is said to
be the characteristic polynomial of ${\rm {\bf A}}$.
\end{definition}
The roots of the characteristic polynomial of a Hermitian matrix
are its real left eigenvalues, which are its right eigenvalues as
well. We can prove the following theorem by analogy to the
commutative case (see, e.g. \cite{la}).
\begin{theorem}\label{theor:char_polin}
If ${\rm {\bf A}} \in {\rm M}\left( {n,{\rm {\mathbb{H}}}}
\right)$ is Hermitian, then $p_{{\rm {\bf A}}}\left( {t} \right) =
t^{n} - d_{1} t^{n - 1} + d_{2} t^{n - 2} - \ldots + \left( { - 1}
\right)^{n}d_{n}$, where $d_{k} $ is the sum of principle minors
of ${\rm {\bf A}}$ of order $rk$, $1 \le k < n$, and $d_{n}=\det
{\rm {\bf A}}$.
\end{theorem}
\section{Determinantal representations of the W-weigh\-ted Drazin inverse for an arbitrary matrix}
Determinantal representations W-weighted Drazin inverse of complex matrices have been received by full-rank factorization in \cite{liu1} and by a limit representation in \cite{liu2}.
For an arbitrary matrix over the field of complex numbers, ${\bf A}\in  {\mathbb{C}}^{m\times n}$, we denote by $\mathcal{R}({\rm {\bf A}})$  the
range of ${\rm {\bf A}}$ and by $\mathcal{N}({\rm {\bf A}})$  the null
space of ${\rm {\bf A}}$.
For an arbitrary matrix over the quaternion skew field, ${\bf A}\in  {\mathbb{H}}^{m\times n}$, we denote by

$
 \mathcal{R}_{r}({\rm {\bf A}})=\{ {\bf y}\in {\mathbb{H}}^{m} : \,\,{\bf y} = {\bf A}{\bf x},\,\,  {\bf x} \in {\mathbb{H}}^{n}\},$ the column right space of ${\bf A}$,

$\mathcal{N}_{r}({\rm {\bf A}})=\{ {\bf y}\in {\mathbb{H}}^{n} : \,\, {\bf A}{\bf x}=0\}$,  the right null space of  ${\bf A}$,

$\mathcal{R}_{l}({\rm {\bf A}})=\{ {\bf y}\in {\mathbb{H}}^{n} : \,\,{\bf y} = {\bf x}{\bf A},\,\,  {\bf x} \in {\mathbb{H}}^{m}\},$ the column left space of ${\bf A}$,

$\mathcal{N}_{r}({\rm {\bf A}})=\{ {\bf y}\in {\mathbb{H}}^{m} : \,\, {\bf x}{\bf A}=0\}$,  the left null space of  ${\bf A}.$

Through the theory of the column-row determinants, a determinantal representation W-weighted Drazin inverse over the quaternion skew-field for the first time has been obtained in \cite{song2} by the following theorem.
\begin{theorem}
Let ${\bf A}\in {\mathbb{H}}^{m\times n}$,  ${\bf W}\in {\mathbb{H}}^{n\times m}$ with $k= {\rm max}\{Ind({\bf A}{\bf W}), Ind({\bf W}{\bf A})\}$ and and $
\rank({\bf A}{\bf W})^{k}=s$. Suppose that ${\bf B}\in {\mathbb{H}}^{n\times (n-s)}_{n-s}$ and ${\bf C}^{*}\in {\mathbb{H}}^{m\times (m-s)}_{m-s}$ are of full ranks and
\begin{gather*}
\mathcal{R}_{r}( {\bf B})=\mathcal{N}_{r}\left(({\bf W}{\bf A})^{k}\right),\,\,\,\,\mathcal{N}_{r}( {\bf C})=\mathcal{R}_{r}\left(({\bf A}{\bf W})^{k}\right),\\
\mathcal{R}_{l}( {\bf C})=\mathcal{N}_{l}\left(({\bf A}{\bf W})^{k}\right),\,\,\,\,\mathcal{N}_{l}( {\bf B})=\mathcal{R}_{l}\left(({\bf W}{\bf A})^{k}\right).
\end{gather*}
Denote
\begin{gather*}
{\bf M}=\begin{bmatrix} {\bf W}{\bf A}{\bf W} & {\bf B} \\ {\bf C} & 0 \end{bmatrix}.
\end{gather*}
Then the W-weighted Drazin inverse ${\rm {\bf A}}_{d,{\bf W}}=(a)_{ij}\in {\mathbb{H}}^{n\times m}$ has the following determinantal representations:
\begin{equation}\label{eq:WDr_so1}
a_{ij}= \frac{\sum^{m+n-s}_{k=1}L_{ki}m^{*}_{kj}}{\det{\bf M}^{*}{\bf M}}, i=\overline{1,m}, j=\overline{1,n},
\end{equation}
or
\begin{equation}\label{eq:WDr_so2}
a_{ij}= \frac{\sum^{m+n-s}_{k=1}m^{*}_{ik}R_{jk}}{\det{\bf M}{\bf M}^{*}}, i=\overline{1,m}, j=\overline{1,n},
\end{equation}
where $L_{ij}$ are the left $(ij)$-th cofactor of  ${\bf M}^{*}{\bf M}$ and $R_{ij}$ are the right $(ij)$-th cofactor of  ${\bf M}{\bf M}^{*}$, respectively, for all $i,j=\overline{1,m+n-s}$.
\end{theorem}
As can be seen,  the auxiliary matrices ${\bf B}$ and ${\bf C}$ have been used in the determinantal representations (\ref{eq:WDr_so1}) and (\ref{eq:WDr_so2}). In this paper we escape it.
Below we give determinantal representations of the W-weighted Drazin inverse of an arbitrary matrix ${\bf A}\in {\mathbb{H}}^{m\times n}$ with with respect to the matrix ${\bf W}\in {\mathbb{H}}^{n\times m}$
 by using the determinantal representations of the  Drazin
inverse,  of the Moore-Penrose inverse, and the limit representation of the W-weighted Drazin inverse in some particular case.
\subsection{Determinantal representations of the W-weighted\\ Drazin inverse by using determinantal representations of the Drazin inverse}
Denote ${\bf W}{\bf A}=:{\bf U}=
({u}_{ij})\in {\mathbb{H}}^{n\times n}$ and ${\bf A}{\bf W}=:{\bf V}=
({v}_{ij})\in {\mathbb{H}}^{m\times m}$.
Due to Theorem \ref{theor:det_rep_draz}, we denote an entry of the Drazin inverse ${\bf U}^D$ by
\begin{equation}
\label{eq:u_cdet_draz} u_{ij} ^{D,1}=
 {\frac{{ \sum\limits_{t = 1}^{n} {u}_{it}^{(k)}   {\sum\limits_{\beta \in J_{r,\,n} {\left\{ {t}
\right\}}} {{\rm{cdet}} _{t} \left( {\left({\bf U}^{ 2k+1} \right)^{*}\left({\bf U}^{ 2k+1} \right)_{. \,t} \left( \hat{{\rm {\bf u}}}_{.\,j}
\right)} \right){\kern 1pt}  _{\beta} ^{\beta} } }
}}{{{\sum\limits_{\beta \in J_{r,\,n}} {{\left| {\left({\bf U}^{ 2k+1} \right)^{*}\left({\bf U}^{ 2k+1} \right){\kern 1pt} _{\beta} ^{\beta}
}  \right|}}} }}}
\end{equation}
or
\begin{equation}
\label{eq:u_rdet_draz} u_{ij} ^{D,2}=
{\frac{\sum\limits_{s = 1}^{n}\left({{\sum\limits_{\alpha \in I_{r,\,n} {\left\{ {s}
\right\}}} {{\rm{rdet}} _{s} \left( {\left( { {\bf U}^{ 2k+1} \left({\bf U}^{ 2k+1} \right)^{*}
} \right)_{\,. s} (\check{{\rm {\bf u}}}_{i\,.})} \right) {\kern 1pt} _{\alpha} ^{\alpha} } }
}\right){u}_{sj}^{(k)}}{{{\sum\limits_{\alpha \in I_{r,\,n}} {{\left| {\left( { {\bf U}^{ 2k+1} \left({\bf U}^{ 2k+1} \right)^{*}
} \right){\kern 1pt} _{\alpha} ^{\alpha}
}  \right|}}} }}}
\end{equation}
where $\hat{{\rm {\bf u}}}_{.s}$ and $\check{{{\rm {\bf u}}}}_{t.}$ the $s$th column of $
({\bf U}^{ 2k+1})^{*} {\bf U}^{k}=:\hat{{\rm {\bf U}}}=
(\hat{u}_{ij})\in {\mathbb{H}}^{n\times n}$
and  the $t$th row of $
{\bf U}^{k}({\bf U}^{ 2k+1})^{*} =:\check{{\rm {\bf U}}}=
(\check{u}_{ij})\in {\mathbb{H}}^{n\times n}$, respectively  for all $s,t=\overline
{1,n}$, $r=\rank{\rm {\bf U}}^{k+1} =
\rank{\rm {\bf U}}^{k}$.
Then by  (\ref{eq:WDr_dr}) we can obtain the following determinantal representations of ${\rm {\bf A}}_{d,{\bf W}}=
({a}_{ij}^{d,{\bf W}})\in {\mathbb{H}}^{m\times n}$,
\begin{equation}
\label{eq:det_rep1_wdraz}
{a}_{ij}^{d,{\bf W}}=\sum\limits_{q = 1}^{n} a_{iq}(u_{qj}^{D})^{(2)}
\end{equation}
where $(u_{qj}^{D})^{(2)}=\sum\limits_{p = 1}^{n}u_{qp}^{D,l}u_{pj}^{D,f}$ for all $l,f=\overline
{1,2}$, and $u_{ij} ^{D,1}$ from (\ref{eq:u_cdet_draz}) and $u_{ij} ^{D,2}$ from (\ref{eq:u_rdet_draz}).
Similarly using ${\bf V}=
({v}_{ij})\in {\mathbb{H}}^{m\times m}$ we have the following determinantal representations of ${\rm {\bf A}}_{d,{\bf W}}$,
\begin{equation}
\label{eq:det_rep2_wdraz}
{a}_{ij}^{d,{\bf W}}=\sum\limits_{q = 1}^{m} (v_{iq}^{D})^{(2)}a_{qj}.
\end{equation}
The first factor is one of the four possible equations $(v_{iq}^{D})^{(2)}=\sum\limits_{p = 1}^{m}v_{ip}^{D,l}v_{pq}^{D,f}$ for all $l,f=\overline
{1,2}$, and an entry of the Drazin inverse ${\bf V}^D$ is denoting by
\begin{equation}
\label{eq:v_cdet_draz} v_{ij} ^{D,1}=
 {\frac{{ \sum\limits_{t = 1}^{m} {v}_{it}^{(k)}   {\sum\limits_{\beta \in J_{r,\,m} {\left\{ {t}
\right\}}} {{\rm{cdet}} _{t} \left( {\left({\bf V}^{ 2k+1} \right)^{*}\left({\bf V}^{ 2k+1} \right)_{. \,t} \left( \hat{{\rm {\bf v}}}_{.\,j}
\right)} \right){\kern 1pt}  _{\beta} ^{\beta} } }
}}{{{\sum\limits_{\beta \in J_{r,\,m}} {{\left| {\left({\bf V}^{ 2k+1} \right)^{*}\left({\bf V}^{ 2k+1} \right){\kern 1pt} _{\beta} ^{\beta}
}  \right|}}} }}}
\end{equation}
or
\begin{equation}
\label{eq:v_rdet_draz} v_{ij} ^{D,2}=
{\frac{\sum\limits_{s = 1}^{m}\left({{\sum\limits_{\alpha \in I_{r,\,m} {\left\{ {s}
\right\}}} {{\rm{rdet}} _{s} \left( {\left( { {\bf V}^{ 2k+1} \left({\bf V}^{ 2k+1} \right)^{*}
} \right)_{\,. s} (\check{{\rm {\bf v}}}_{i\,.})} \right) {\kern 1pt} _{\alpha} ^{\alpha} } }
}\right){v}_{sj}^{(k)}}{{{\sum\limits_{\alpha \in I_{r,\,m}} {{\left| {\left( { {\bf V}^{ 2k+1} \left({\bf V}^{ 2k+1} \right)^{*}
} \right){\kern 1pt} _{\alpha} ^{\alpha}
}  \right|}}} }}},
\end{equation}
where $\hat{{\rm {\bf v}}}_{.s}$ and $\check{{{\rm {\bf v}}}}_{t.}$ the $s$th column of $
({\bf V}^{ 2k+1})^{*} {\bf V}^{k}=:\hat{{\rm {\bf V}}}=
(\hat{v}_{ij})\in {\mathbb{H}}^{m\times m}$
and  the $t$th row of $
{\bf V}^{k}({\bf V}^{ 2k+1})^{*} =:\check{{\rm {\bf V}}}=
(\check{v}_{ij})\in {\mathbb{H}}^{m\times m}$, respectively  for all $s,t=\overline
{1,m}$, $r=\rank{\rm {\bf V}}^{k+1} =
\rank{\rm {\bf V}}^{k}$.
\subsection{Determinantal representations of the W-weighted\\ Drazin
inverse by using determinantal representations of the Moore-Penrose inverse}
Consider the general algebraic structures (GAS) of the matrices  ${\bf A}\in {\mathbb{H}}^{m\times n}$, ${\bf W}\in {\mathbb{H}}^{n\times m}$, ${\bf A}^{+}\in {\mathbb{H}}^{n\times m}$, ${\bf W}^{+}\in {\mathbb{H}}^{m\times n}$ and ${\rm {\bf A}}_{d,{\bf W}}\in {\mathbb{H}}^{m\times n}$ with $k= {\rm max}\{Ind({\bf A}{\bf W}), Ind({\bf W}{\bf A})\}$ (e.g.,\cite{wei1,wei2,wei4,ze}).
Let exist ${\bf L}\in {\mathbb{H}}^{m\times m}$ and ${\bf Q}\in {\mathbb{H}}^{n\times n}$ such that
\[
{\bf A}={\bf L} \left[
                        \begin{array}{cc}
                         {\bf A}_{11} & {\bf 0} \\
                          {\bf 0} &  {\bf A}_{22} \\
                        \end{array}
                      \right]{\bf Q}^{-1},\,\,\,\,{\bf W}={\bf Q} \left[
                        \begin{array}{cc}
                         {\bf W}_{11} & {\bf 0} \\
                          {\bf 0} &  {\bf W}_{22} \\
                        \end{array}
                      \right]{\bf L}^{-1}.\]
Then \[
{\bf A}^{+}= {\bf Q}\left[
                        \begin{array}{cc}
                         {\bf A}_{11}^{-1} & {\bf 0} \\
                          {\bf 0} & {\bf 0} \\
                        \end{array}
                      \right]{\bf L}^{-1},\,\,\,\,{\bf W}^{+}={\bf L} \left[
                        \begin{array}{cc}
                         {\bf W}_{11}^{-1} & {\bf 0} \\
                          {\bf 0} &  {\bf 0} \\
                        \end{array}
                      \right]{\bf Q}^{-1},\]
 \[
{\bf A}_{d,{\bf W}}= {\bf L}\left[
                        \begin{array}{cc}
                        ( {\bf W}_{11} {\bf A}_{11} {\bf W}_{11})^{-1} & {\bf 0} \\
                          {\bf 0} & {\bf 0} \\
                        \end{array}
                      \right]{\bf Q}^{-1},\]
where ${\bf L}, {\bf Q}, {\bf A}_{11}, {\bf W}_{11}$ are non-singular matrices, and $ {\bf A}_{22}{\bf W}_{22}, {\bf W}_{22}{\bf A}_{22}$ are nilpotent matrices.
The follow theorem due to \cite{ze} can be expanded to  ${\mathbb{H}}$.
\begin{theorem}\label{theor:al1_wdraz_det}Let ${\bf A}\in {\mathbb{H}}^{m\times n}$ and ${\bf W}\in {\mathbb{H}}^{n\times m}$  such that ${\bf A}_{22}{\bf W}_{22}$ and ${\bf W}_{22}{\bf A}_{22}$ are
nilpotent matrices of index $k$ in GAS form. Then  the weighted Drazin inverse of ${\bf A}$ with respect to ${\bf W}$ can be written as matrix expression involving the Moore-Penrose inverse,
\begin{equation}
\label{eq:al1_wdraz_det}
{\bf A}_{d,{\bf W}}= \left\{ ({\bf A}{\bf W})^{k}\left[({\bf A}{\bf W})^{2k+1}\right]^{+}({\bf A}{\bf W})^{k}\right\}{\bf W}^{+},
\end{equation}
where $k= {\rm max}\{Ind({\bf A}{\bf W}), Ind({\bf W}{\bf A})\}$.
\end{theorem}
Similarly can be obtained the following theorem.
\begin{theorem}\label{theor:al2_wdraz_det}Let ${\bf A}\in {\mathbb{H}}^{m\times n}$ and ${\bf W}\in {\mathbb{H}}^{n\times m}$  such that ${\bf A}_{22}{\bf W}_{22}$ and ${\bf W}_{22}{\bf A}_{22}$ are
nilpotent matrices of index $k$ in GAS form. Then  the W-weighted Drazin inverse of ${\bf A}$ with respect to ${\bf W}$ can be written as the following matrix expression,
\begin{equation}
\label{eq:al2_wdraz_det}
{\bf A}_{d,{\bf W}}= {\bf W}^{+}\left\{ ({\bf W}{\bf A})^{k}\left[({\bf W}{\bf A})^{2k+1}\right]^{+}({\bf W}{\bf A})^{k}\right\},
\end{equation}
where $k= {\rm max}\{Ind({\bf A}{\bf W}), Ind({\bf W}{\bf A})\}$.
\end{theorem}
{\bf Proof}. Since ${\bf W}_{22}{\bf A}_{22}$ is a nilpotent matrix of index $k$, then due to the GAS of ${\bf A}$, ${\bf W}$ and their generalized inverses we have the following Jordan canonical forms,
\[
{\bf W}{\bf A}={\bf Q} \left[
                        \begin{array}{cc}
                         {\bf W}_{11}{\bf A}_{11} & {\bf 0} \\
                          {\bf 0} &  {\bf W}_{22}{\bf A}_{22} \\
                        \end{array}
                      \right]{\bf Q}^{-1},\,\,\,\,({\bf W}{\bf A})^{k}={\bf Q} \left[
                        \begin{array}{cc}
                        ({\bf W}_{11}{\bf A}_{11})^{k} & {\bf 0} \\
                          {\bf 0} &  {\bf 0} \\
                        \end{array}
                      \right]{\bf Q}^{-1},\]
\[\left[({\bf W}{\bf A})^{2k+1}\right]^{+}={\bf Q} \left[
                        \begin{array}{cc}
                        ({\bf W}_{11}{\bf A}_{11})^{-2k-1} & {\bf 0} \\
                          {\bf 0} &  {\bf 0} \\
                        \end{array}
                      \right]{\bf Q}^{-1}.\]
Simple computing of ${\bf W}^{+}\left\{ ({\bf W}{\bf A})^{k}\left[({\bf W}{\bf A})^{2k+1}\right]^{+}({\bf W}{\bf A})^{k}\right\}$ proves the theorem,
\begin{multline*}
{\bf W}^{+}\left\{ ({\bf W}{\bf A})^{k}\left[({\bf W}{\bf A})^{2k+1}\right]^{+}({\bf W}{\bf A})^{k}\right\}=\\
{\bf L} \begin{bmatrix}
                         {\bf W}_{11}^{-1} & {\bf 0} \\
                          {\bf 0} &  {\bf 0} \\
                        \end{bmatrix}
                      \begin{bmatrix}
                        ({\bf W}_{11}{\bf A}_{11})^{k} & {\bf 0} \\
                          {\bf 0} &  {\bf 0} \\
                         \end{bmatrix}
                      \begin{bmatrix}
                        ({\bf W}_{11}{\bf A}_{11})^{-2k-1} & {\bf 0} \\
                          {\bf 0} &  {\bf 0} \\
                         \end{bmatrix}
                      \begin{bmatrix}
                        ({\bf W}_{11}{\bf A}_{11})^{k} & {\bf 0} \\
                          {\bf 0} &  {\bf 0} \\
                         \end{bmatrix}{\bf Q}^{-1}=\\
                         {\bf L} \begin{bmatrix}
                         {\bf W}_{11}^{-1}({\bf W}_{11}{\bf A}_{11})^{k}({\bf W}_{11}{\bf A}_{11})^{-2k-1}({\bf W}_{11}{\bf A}_{11})^{k} & {\bf 0} \\
                          {\bf 0} &  {\bf 0} \\
                        \end{bmatrix}{\bf Q}^{-1}=\\
                        {\bf L} \begin{bmatrix}
                         {\bf W}_{11}^{-1}({\bf W}_{11}{\bf A}_{11})^{-1} & {\bf 0} \\
                          {\bf 0} &  {\bf 0} \\
                        \end{bmatrix}{\bf Q}^{-1}=\\
                        {\bf L} \begin{bmatrix}
                         ({\bf W}_{11}{\bf A}_{11}{\bf W}_{11})^{-1} & {\bf 0} \\
                          {\bf 0} &  {\bf 0} \\
                        \end{bmatrix}{\bf Q}^{-1}={\bf A}_{d,{\bf W}}.\blacksquare
\end{multline*}
Using (\ref{eq:al1_wdraz_det}),  an entry ${a}_{ij}^{d,{\bf W}}$ of the W-weighted Drazin inverse ${\bf A}_{d,{\bf W}}$ can be obtained as follows
\begin{equation}
\label{eq:det_wdraz} {a}_{ij}^{d,{\bf W}}= {\sum\limits_{s = 1}^{m}
{{\sum\limits_{t = 1}^{m}\sum\limits_{l = 1}^{m} {{v}_{is}^{(k)} \left(v_{st}^{(2k+1)}  \right)^{+}} } }
{v}_{tl}^{(k)}}{w}_{lj}^{+}
\end{equation}
 for all $i=\overline {1,m}$, $j =
\overline {1,n}$.
Denote  by $\check{{{\rm {\bf w}}}}_{t.}$ the $t$th row of $
{\bf V}^{k}{\bf W}^{*} =:\check{{\rm {\bf W}}}=
(\check{w}_{ij})\in {\mathbb{H}}^{m\times n}$ for all $t=\overline
{1,m}$. It follows from $\sum\limits_{l}v_{tl}^{(k)}{\bf w}_{l\,.}^{*}=\check{{\rm {\bf w}}}_{t\,.}$ and (\ref{eq:det_repr_AA*}) that
\begin{multline}\label{eq: tow1}
\sum\limits_{l = 1}^{m}{v}_{tl}^{(k)}{w}_{lj}^{+}
=
\sum\limits_{l =
1}^{m}{v}_{tl}^{(k)}\cdot{\frac{{{\sum\limits_{\alpha \in I_{r_{1},\,n} {\left\{ {j}
\right\}}} {{\rm{rdet}} _{j} \left( {\bf W} {\bf W}^{*}
\right)_{j. } \left( { {\bf  w}_{l.}^{*}} \right) {\kern 1pt} _{\alpha} ^{\alpha} } }
}}{{{\sum\limits_{\alpha \in I_{r_{1},\,n}} {{\left| {\left( { {\bf W} {\bf W}^{*}
} \right){\kern 1pt} _{\alpha} ^{\alpha}
}  \right|}}} }}}
=\\
{\frac{{{\sum\limits_{\alpha \in I_{r_{1},\,n} {\left\{ {j}
\right\}}} {{\rm{rdet}} _{j} \left( {\left(  {\bf W} {\bf W}^{*} \right)_{j\,. } (\check{{\rm {\bf w}}}_{t.})} \right) {\kern 1pt} _{\alpha} ^{\alpha} } }
}}{{{\sum\limits_{\alpha \in I_{r_{1},\,n}} {{\left| {\left( {\bf W} {\bf W}^{*}\right){\kern 1pt} _{\alpha} ^{\alpha}
}  \right|}}} }}},
\end{multline}
where $r_{1}=\rank{\rm {\bf W}}$.
Similarly, denote  by $\check{{{\rm {\bf v}}}}_{i.}$ the $t$th row of $
{\bf V}^{k}({\bf V}^{ 2k+1})^{*} =:\check{{\rm {\bf V}}}=
(\check{v}_{ij})\in {\mathbb{H}}^{m\times m}$ for all $t=\overline
{1,m}$. It follows from ${\sum\limits_{s} v_{is}^{(k)}{\left( {\rm {\bf v}}_{s\,.}^{
 (2k+1)} \right)^{*}} }=\check{{\rm {\bf v}}}_{i\,.}$ and (\ref{eq:det_repr_AA*}) that
\begin{multline}\label{eq: tov1}
\sum\limits_{s = 1}^{m}{v}_{is}^{(k)} \left(v_{st}^{(2k+1)}  \right)^{+}
=\sum\limits_{s =
1}^{m}{v}_{is}^{(k)}\cdot{\frac{{{\sum\limits_{\alpha \in I_{r,\,m} {\left\{ {t}
\right\}}} {{\rm{rdet}} _{t} \left( {\left( { {\bf V}^{ 2k+1} \left({\bf V}^{ 2k+1} \right)^{*}
} \right)_{t. } \left( {{\rm {\bf v}}_{s.}^{(2k+1)} }
\right)^{*}} \right) {\kern 1pt} _{\alpha} ^{\alpha} } }
}}{{{\sum\limits_{\alpha \in I_{r,\,m}} {{\left| {\left( { {\bf V}^{ 2k+1} \left({\bf V}^{ 2k+1} \right)^{*}
} \right){\kern 1pt} _{\alpha} ^{\alpha}
}  \right|}}} }}}
=\\
{\frac{{{\sum\limits_{\alpha \in I_{r,\,m} {\left\{ {t}
\right\}}} {{\rm{rdet}} _{t} \left( {\left( { {\bf V}^{ 2k+1} \left({\bf V}^{ 2k+1} \right)^{*}
} \right)_{t. } (\check{{\rm {\bf v}}}_{i\,.})} \right) {\kern 1pt} _{\alpha} ^{\alpha} } }
}}{{{\sum\limits_{\alpha \in I_{r,\,m}} {{\left| {\left( { {\bf V}^{ 2k+1} \left({\bf V}^{ 2k+1} \right)^{*}
} \right){\kern 1pt} _{\alpha} ^{\alpha}
}  \right|}}} }}},
\end{multline}
where $r=\rank{\rm {\bf W}}^{k+1} =
\rank{\rm {\bf W}}^{k}$.
Using (\ref{eq: tow1}) and (\ref{eq: tov1}) in (\ref{eq:det_wdraz})  we obtain the following determinantal representation of ${\bf A}_{d,{\bf W}}$ ,
\begin{multline}
\label{eq:det_repr_v wdraz} {a}_{ij}^{d,{\bf W}}=\\
\frac{\sum\limits_{t = 1}^{m}{\sum\limits_{\alpha \in I_{r,\,m} {\left\{ {t}
\right\}}} {{\rm{rdet}} _{t} \left( {\left( { {\bf V}^{ 2k+1} \left({\bf V}^{ 2k+1} \right)^{*}
} \right)_{t. } (\check{{\rm {\bf v}}}_{i\,.})} \right) {\kern 1pt} _{\alpha} ^{\alpha} } }
{\sum\limits_{\alpha \in I_{r_{1},\,n} {\left\{ {j}
\right\}}} {{\rm{rdet}} _{j} \left( {\left(  {\bf W} {\bf W}^{*} \right)_{j\,. } (\check{{\rm {\bf w}}}_{t.})} \right) {\kern 1pt} _{\alpha} ^{\alpha} } }}{{\sum\limits_{\alpha \in I_{r,\,m}} {{\left| {\left( { {\bf V}^{ 2k+1} \left({\bf V}^{ 2k+1} \right)^{*}
} \right){\kern 1pt} _{\alpha} ^{\alpha}
}  \right|}}}{\sum\limits_{\alpha \in I_{r_{1},\,n}} {{\left| {\left( {\bf W} {\bf W}^{*}\right){\kern 1pt} _{\alpha} ^{\alpha}
}  \right|}}}}
\end{multline}
Thus we have proved the following theorem.
\begin{theorem}\label{theor:det_rep_vwdraz} Let ${\bf A}\in {\mathbb{H}}^{m\times n}$ and ${\bf W}\in {\mathbb{H}}_{r_{1}}^{n\times m}$ with $k= {\rm max}\{Ind({\bf A}{\bf W}), Ind({\bf W}{\bf A})\}$ and $r=\rank({\bf A}{\bf W})^{k+1} =
\rank({\bf A}{\bf W})^{k}$.  Then  the W-weighted Drazin inverse of ${\bf A}$ with respect to ${\bf W}$ possesses the determinantal representation  (\ref{eq:det_repr_v wdraz}), where ${\bf V}={\bf A}{\bf W}$,
$\check{{\rm {\bf V}}}={\bf V}^{k}({\bf V}^{ 2k+1})^{*}$, and $
\check{{\rm {\bf W}}}={\bf V}^{k}{\bf W}^{*}$.
\end{theorem}
Similarly we have the following theorem.
\begin{theorem}\label{theor:det_rep_uwdraz} Let ${\bf A}\in {\mathbb{H}}^{m\times n}$ and ${\bf W}\in {\mathbb{H}}_{r_{1}}^{n\times m}$ with $k= {\rm max}\{Ind({\bf A}{\bf W}), Ind({\bf W}{\bf A})\}$ and $r=\rank({\bf W}{\bf A})^{k+1} =
\rank({\bf W}{\bf A})^{k}$.  Then  the W-weighted Drazin inverse of ${\bf A}$ with respect to ${\bf W}$ possesses the following determinantal representation,
\begin{multline}
\label{eq:det_repr_u wdraz} {a}_{ij}^{d,{\bf W}}=\\
\frac{\sum\limits_{t = 1}^{n}
{{\sum\limits_{\beta \in J_{r_{1},\,m} {\left\{ {i}
\right\}}} {{\rm{cdet}} _{i} \left( {\left(  {\bf W}^{*} {\bf W} \right)_{. i} (\hat{{\rm {\bf w}}}_{.t})} \right) {\kern 1pt} _{\beta} ^{\beta} } }
}
{{\sum\limits_{\beta \in J_{r,\,n} {\left\{ {t}
\right\}}} {{\rm{cdet}} _{t} \left( {\left( { \left({\bf U}^{ 2k+1} \right)^{*}{\bf U}^{ 2k+1}
} \right)_{.t } (\hat{{\rm {\bf u}}}_{.j})} \right) {\kern 1pt} _{\beta} ^{\beta} } }
}}
{{\sum\limits_{\beta \in J_{r_{1},\,m}} {{\left| {\left( {\bf W}^{*} {\bf W}\right){\kern 1pt} _{\beta} ^{\beta}
}  \right|}}}{\sum\limits_{\beta \in J_{r,\,
n}} {{\left| {\left( { \left({\bf U}^{ 2k+1} \right)^{*}{\bf U}^{ 2k+1}
} \right){\kern 1pt} _{\beta} ^{\beta}
}  \right|}}}}
\end{multline}
where ${\bf U}={\bf W}{\bf A}$,
$\hat{{\rm {\bf U}}}=({\bf U}^{ 2k+1})^{*}{\bf U}^{k}$, and $
\hat{{\rm {\bf W}}}={\bf W}^{*}{\bf U}^{k}$.
\end{theorem}
{\bf Proof}.
Using (\ref{eq:al2_wdraz_det}),  an entry ${a}_{ij}^{d,{\bf W}}$ of the W-weighted Drazin inverse ${\bf A}_{d,{\bf W}}$ can be obtained as follows
\begin{equation}
\label{eq:det_wdraz2} {a}_{ij}^{d,{\bf W}}= {\sum\limits_{s = 1}^{n}
{{\sum\limits_{t = 1}^{n}\sum\limits_{l = 1}^{
n}{w}_{is}^{+} {{u}_{st}^{(k)} \left(u_{tl}^{(2k+1)}  \right)^{+}} } }
{u}_{lj}^{(k)}}
\end{equation}
 for all $i=\overline {1,m}$, $j =
\overline {1,n}$.
Denote  by $\hat{{{\rm {\bf w}}}}_{.t}$ the $t$th column of $
{\bf W}^{*}{\bf U}^{k} =:\hat{{\rm {\bf W}}}=
(\hat{w}_{ij})\in {\mathbb{H}}^{m\times n}$ for all $t=\overline
{1,n}$.  It follows from $\sum\limits_{t}{\bf w}_{.s}^{*}u_{st}^{(k)}=\hat{ {\bf w}}_{.t}$ and (\ref{eq:det_repr_A*A}) that
\begin{multline}\label{eq: tow2}
\sum\limits_{s = 1}^{n}{w}_{is}^{+}{u}_{st}^{(k)}
=
\sum\limits_{s =
1}^{n}{\frac{{{\sum\limits_{\beta \in J_{r_{1},\,m} {\left\{ {i}
\right\}}} {{\rm{cdet}} _{i} \left( {\bf W}^{*}{\bf W}
\right)_{.i } \left( { {\bf  w}_{.s}^{*}} \right) {\kern 1pt} _{\beta} ^{\beta} } }
}}{{{\sum\limits_{\beta \in J_{r_{1},\,m}} {{\left| {\left( { {\bf W}^{*}{\bf W}
} \right){\kern 1pt} _{\beta} ^{\beta}
}  \right|}}} }}}\cdot{u}_{st}^{(k)}
=\\
{\frac{{{\sum\limits_{\beta \in J_{r_{1},\,m} {\left\{ {i}
\right\}}} {{\rm{cdet}} _{i} \left( {\left(  {\bf W}^{*} {\bf W} \right)_{. i} (\hat{{\rm {\bf w}}}_{.t})} \right) {\kern 1pt} _{\beta} ^{\beta} } }
}}{{{\sum\limits_{\beta \in J_{r_{1},\,m}} {{\left| {\left( {\bf W}^{*} {\bf W}\right){\kern 1pt} _{\beta} ^{\beta}
}  \right|}}} }}},
\end{multline}
where $r_{1}=\rank{\rm {\bf W}}$.
Similarly, denote  by $\hat{{{\rm {\bf u}}}}_{.j}$ the $j$th column of $
({\bf U}^{ 2k+1})^{*}{\bf U}^{k} =:\hat{{\rm {\bf U}}}=
(\hat{u}_{ij})\in {\mathbb{H}}^{n\times n}$ for all $j=\overline
{1,n}$. It follows from ${\sum\limits_{l}{\left( {\rm {\bf u}}_{.l}^{
 (2k+1)} \right)^{*}} } u_{lj}^{(k)}=\hat{{\rm {\bf u}}}_{.j}$ and (\ref{eq:det_repr_A*A}) that
 \begin{multline}\label{eq: tou2}
\sum\limits_{l = 1}^{n} \left(u_{tl}^{(2k+1)}  \right)^{+}{u}_{lj}^{(k)}
=\\\sum\limits_{l =
1}^{n}{\frac{{{\sum\limits_{\beta \in J_{r,\,n} {\left\{ {t}
\right\}}} {{\rm{cdet}} _{t} \left( {\left( {  \left({\bf U}^{ 2k+1} \right)^{*}{\bf U}^{ 2k+1}
} \right)_{.t } \left( {{\rm {\bf u}}_{.l}^{(2k+1)} }
\right)^{*}} \right) {\kern 1pt} _{\beta} ^{\beta} } }
}}{{{\sum\limits_{\beta \in J_{r,\,n}} {{\left| {\left( { \left({\bf U}^{ 2k+1} \right)^{*}{\bf U}^{ 2k+1}
} \right){\kern 1pt} _{\beta} ^{\beta}
}  \right|}}} }}}\cdot{u}_{lj}^{(k)}
=\\
{\frac{{{\sum\limits_{\beta \in J_{r,\,n} {\left\{ {t}
\right\}}} {{\rm{cdet}} _{t} \left( {\left( { \left({\bf U}^{ 2k+1} \right)^{*}{\bf U}^{ 2k+1}
} \right)_{.t } (\hat{{\rm {\bf u}}}_{.j})} \right) {\kern 1pt} _{\beta} ^{\beta} } }
}}{{{\sum\limits_{\beta \in J_{r,\,n}} {{\left| {\left( { \left({\bf U}^{ 2k+1} \right)^{*}{\bf U}^{ 2k+1}
} \right){\kern 1pt} _{\beta} ^{\beta}
}  \right|}}} }}},
\end{multline}
where $r=\rank({\bf A}{\bf W})^{k+1} =
\rank({\bf A}{\bf W})^{k}$. Using the equations (\ref{eq: tou2}) and (\ref{eq: tow2}) in
 (\ref{eq:det_wdraz2}), we obtain (\ref{eq:det_repr_u wdraz}).
$\blacksquare$
\subsection{Determinantal representations of the W-weighted\\ Drazin
inverse in some special case}
In this subsection we consider  the determinantal representation of the
W-weighted Drazin inverse of  ${\bf A}\in {\mathbb{H}}^{m\times n}$  with
respect to  ${\bf W}\in {\mathbb{H}}^{n\times m}$ in a special case, when
 ${\bf A}{\bf W}={\bf V}=\left(v_{ij}
\right)  \in {\mathbb{H}}^{m\times m}$ and ${\bf W}{\bf A}={\bf
U}=\left(u_{ij} \right)\in {\mathbb{H}}^{n\times n}$ are Hermitian. Then for
the determinantal representations
 their Drazin inverse we can use (\ref{eq:dr_rep_cdet}) and (\ref{eq:dr_rep_rdet}).
 For Hermitian matrix, we apply the method, which consists of theorem on the limit representation of the Drazin inverse, lemmas on rank of matrices and on characteristic polynomial. This method was used at first in \cite{ky8} and then in \cite{ky3,liu2,ky9,ky10}.
By analogy to the complex case \cite{car} we have the following  limit representations of the W-weighted
Drazin inverse,
\begin{equation}\label{eq:lim_repr_AW}
{\rm {\bf A}}_{d,W} = {\mathop {\lim} \limits_{\lambda \to 0}}
\left( {\lambda {\rm {\bf I}}_m + ({\rm {\bf A}}{\bf W})^{k + 2}} \right)^{
- 1}({\rm {\bf A}}{\bf W})^{k}{\bf A}
\end{equation}
and
\begin{equation}\label{eq:lim_repr_WA}
{\rm {\bf A}}_{d,W} = {\mathop {\lim} \limits_{\lambda \to 0}}{\bf A}({\bf W}{\rm {\bf A}})^{k}
\left( {\lambda {\rm {\bf I}}_n + ({\bf W}{\rm {\bf A}})^{k + 2}} \right)^{
- 1}
\end{equation}
where  $\lambda \in {\mathbb R} _{ +}  $, and ${\mathbb R} _{ +} $
is a set of the real positive numbers.

 Denote by ${\rm {\bf v}}_{.j}^{(k)} $ and ${\rm {\bf
v}}_{i.}^{(k)} $ the $j$th column  and the $i$th row of  ${\rm
{\bf V}}^{k} $ respectively. Denote by $ \bar{{\bf V}}^{k}:=({\bf A}{\bf W})^{k}{\bf A}\in {\mathbb{H}}^{m\times n}$ and $ \bar{{\bf W}}={\bf W}{\bf A}{\bf W}\in {\mathbb{H}}^{n\times m}$.
\begin{lemma} \label{lem:rank_col} If ${\bf A}{\bf W}={\bf
V}=\left(v_{ij} \right)\in {\mathbb{H}}^{m\times m}$ with $ Ind{\kern 1pt} {\rm {\bf V}}=k$, then
\begin{equation}\label{eq:rank_col}
 \rank\,\left( {{\rm {\bf V}}^{ k+2} }
\right)_{.\,i} \left( { \bar{{\bf v}}_{.j}^{ (k)} }  \right) \le
\rank\,\left( {{\rm {\bf V}}^{k+2} } \right).
\end{equation}
\end{lemma}
{\bf Proof}. We have ${\rm {\bf V}}^{
 k+2}=\bar{{\bf V}}^{ k}\bar{{\bf W}}$.
Let ${\rm {\bf P}}_{i\,s} \left( {-\bar{w}_{j\,s}}  \right)\in {\mathbb
H}^{m\times m} $, $(s \ne i )$, be a matrix with $-\bar{w}_{j\,s} $ in
the $(i, s)$ entry, 1 in all diagonal entries, and 0 in others.
The  elementary matrix ${\rm {\bf P}}_{i\,s} \left( {-\bar{w}_{j\,s}}  \right)$, $(s \ne i )$, is a matrix of an elementary transformation. It follows that
\[
\left( {\rm {\bf V}}^{ k+2} \right)_{.\,i} \left(
{ \bar{{\bf v}}_{.\,j}^{ (k)} }  \right) \cdot {\prod\limits_{s \ne
i} {{\rm {\bf P}}_{i\,s} \left( {-\bar{w}_{j\,s}}  \right) =
{\mathop
{\left( {{\begin{array}{*{20}c}
 {{\sum\limits_{s \ne j} {\bar{v}_{1s}^{ (k)}  \bar{w}_{s1}} } } \hfill & {\ldots}  \hfill
& {\bar{v}_{1j}^{ (k)} }  \hfill & {\ldots}  \hfill & {{\sum\limits_{s
\ne
j } {\bar{v}_{1s}^{ (k)}  \bar{w}_{sm}}}} \hfill \\
 {\ldots}  \hfill & {\ldots}  \hfill & {\ldots}  \hfill & {\ldots}  \hfill &
{\ldots}  \hfill \\
 {{\sum\limits_{s \ne j} {\bar{v}_{ms}^{(k)}  \bar{w}_{s1}} } } \hfill & {\ldots}  \hfill
& {\bar{v}_{mj}^{(k)} }  \hfill & {\ldots}  \hfill & {{\sum\limits_{s
\ne
j } {\bar{v}_{ms}^{(k)}  \bar{w}_{sm}}}} \hfill \\
\end{array}} }
\right)}\limits_{i-th}}}}.
\]
 We have the next factorization of the obtained matrix.
\[
{\mathop
{\left( {{\begin{array}{*{20}c}
 {{\sum\limits_{s \ne j} {\bar{v}_{1s}^{ (k)}  \bar{w}_{s1}} } } \hfill & {\ldots}  \hfill
& {\bar{v}_{1j}^{ (k)} }  \hfill & {\ldots}  \hfill & {{\sum\limits_{s
\ne
j } {\bar{v}_{1s}^{ (k)}  \bar{w}_{sm}}}} \hfill \\
 {\ldots}  \hfill & {\ldots}  \hfill & {\ldots}  \hfill & {\ldots}  \hfill &
{\ldots}  \hfill \\
 {{\sum\limits_{s \ne j} {\bar{v}_{ms}^{(k)}  \bar{w}_{s1}} } } \hfill & {\ldots}  \hfill
& {\bar{v}_{mj}^{(k)} }  \hfill & {\ldots}  \hfill & {{\sum\limits_{s
\ne
j } {\bar{v}_{ms}^{(k)}  \bar{w}_{sm}}}} \hfill \\
\end{array}} }
\right)}\limits_{i-th}} =
\]
\[
 = \left( {\begin{array}{*{20}c}
 {\bar{v}_{11}^{ (k)} }  \hfill & {\bar{v}_{12}^{ (k)} }  \hfill & {\ldots}  \hfill &
{\bar{v}_{1n}^{ (k)} }  \hfill \\
 {\bar{v}_{21}^{ (k)} }  \hfill & {\bar{v}_{22}^{ (k)} }  \hfill & {\ldots}  \hfill &
{\bar{v}_{2n}^{ (k)} }  \hfill \\
 {\ldots}  \hfill & {\ldots}  \hfill & {\ldots}  \hfill & {\ldots}  \hfill
\\
 {\bar{v}_{m1}^{ (k)} }  \hfill & {\bar{v}_{m2}^{ (k)} }  \hfill & {\ldots}  \hfill &
{\bar{v}_{mn}^{ (k)} }  \hfill \\
\end{array}}  \right)
{\mathop {\left( {{\begin{array}{*{20}c}
 {\bar{w}_{11}}  \hfill & {\ldots}  \hfill & {0} \hfill & {\ldots}  \hfill &
{\bar{w}_{1m}}  \hfill \\
 {\ldots}  \hfill & {\ldots}  \hfill & {\ldots}  \hfill & {\ldots}  \hfill &
{\ldots}  \hfill \\
 {0} \hfill & {\ldots}  \hfill & {1} \hfill & {\ldots}  \hfill & {0} \hfill
\\
 {\ldots}  \hfill & {\ldots}  \hfill & {\ldots}  \hfill & {\ldots}  \hfill &
{\ldots}  \hfill \\
 {\bar{w}_{n1}}  \hfill & {\ldots}  \hfill & {0} \hfill & {\ldots}  \hfill &
{\bar{w}_{nm}}  \hfill \\
\end{array}} } \right)}\limits_{i-th}} j-th.
\]
 Denote ${\rm {\bf \tilde {W}}}: = {\mathop
{\left( {{\begin{array}{*{20}c}
 {\bar{w}_{11}}  \hfill & {\ldots}  \hfill & {0} \hfill & {\ldots}  \hfill &
{\bar{w}_{1m}}  \hfill \\
 {\ldots}  \hfill & {\ldots}  \hfill & {\ldots}  \hfill & {\ldots}  \hfill &
{\ldots}  \hfill \\
 {0} \hfill & {\ldots}  \hfill & {1} \hfill & {\ldots}  \hfill & {0} \hfill
\\
 {\ldots}  \hfill & {\ldots}  \hfill & {\ldots}  \hfill & {\ldots}  \hfill &
{\ldots}  \hfill \\
 {\bar{w}_{n1}}  \hfill & {\ldots}  \hfill & {0} \hfill & {\ldots}  \hfill &
{\bar{w}_{nm}}  \hfill \\
\end{array}} } \right)}\limits_{i-th}} j-th$. The matrix
${\rm {\bf \tilde {W}}}$ is obtained from $ \bar{{\bf W}}={\bf W}{\bf A}{\bf W}$ by
replacing all entries of the $j$th row  and the $i$th column with
zeroes except for 1 in the $(i, j)$ entry. Since elementary
transformations of a matrix do not change a rank, then $\rank{\rm
{\bf V}}^{ k+2} _{.\,i} \left( { \bar{{\bf v}}_{.j}^{ (k)} }
\right) \le \min {\left\{ {\rank \bar{{\bf V}}^{ k },\rank{\rm {\bf
\tilde {W}}}} \right\}}$. It is obvious that
  \[\begin{array}{c}
      \rank \bar{{\bf V}}^{ k }=\rank \,({\bf A}{\bf W})^{k}{\bf A}\geq  \rank\, ({\bf A}{\bf W})^{k+2}, \\
       \rank {\bf \tilde {W}}\geq \rank \,{\bf W}{\bf A}{\bf W}\geq  \rank\, ({\bf A}{\bf W})^{k+2}.
    \end{array}
     \]
 From this the inequality (\ref{eq:rank_col})
follows immediately.
$\blacksquare$

The next lemma is  proved similarly.
\begin{lemma} If ${\bf W}{\bf A}={\bf
U}=\left(u_{ij} \right)\in {\mathbb{H}}^{n\times n}$ with $ Ind{\kern 1pt} {\rm {\bf U}}=k$, then \[ \rank\left( {{\rm
{\bf U}}^{ k+2} } \right)_{i\,.} \left( { \bar{{\bf u}}_{j\,.}^{
(k)} }  \right) \le \rank\left( {{\rm {\bf U}}^{k+2} } \right),\]
where $ \bar{{\bf U}}^{k}:={\bf A}({\bf W}{\bf A})^{k}\in {\mathbb{H}}^{m\times n}$
\end{lemma}
Analogues of the characteristic polynomial are considered in the following two lemmas.
\begin{lemma}\label{lem:char_cdet}
If ${\bf A}{\bf W}={\bf
V}=\left(v_{ij} \right)\in {\mathbb{H}}^{m\times m}$ is Hermitian with $ Ind{\kern 1pt} {\rm {\bf V}}=k$ and $t \in \mathbb{R}$, then
\begin{equation}
\label{eq:char_cdet} {\rm{cdet}} _{i} \left(\lambda {\rm {\bf I}}_m + {{\rm {\bf V}}^{ k+2} }
\right)_{.\,i} \left( { \bar{{\bf v}}_{.j}^{ (k)} }  \right) = c_{1}^{\left( {ij} \right)} \lambda^{n - 1}
+ c_{2}^{\left( {ij} \right)} \lambda^{n - 2} + \ldots + c_{n}^{\left(
{ij} \right)},
\end{equation}
\noindent where $c_{n}^{\left( {ij} \right)} = {\rm{cdet}} _{i}
\left( {{\rm {\bf V}}^{ k+2} }
\right)_{.\,i} \left( { \bar{{\bf v}}_{.j}^{ (k)} }  \right)$  and $c_{s}^{\left( {ij} \right)} =
{\sum\limits_{\beta \in J_{s,\,n} {\left\{ {i} \right\}}}
{{\rm{cdet}} _{i} \left( {\left( {{\rm {\bf V}}^{ k+2} }
\right)_{.\,i} \left( { \bar{{\bf v}}_{.j}^{ (k)} }  \right)}
\right){\kern 1pt}  _{\beta} ^{\beta} } }$ for all $s = \overline
{1,n - 1} $, $i,j = \overline {1,n}$.
\end{lemma}
{\bf Proof}.
Consider the Hermitian matrix $\left( {t{\rm {\bf I}} + {\rm {\bf
V}}^{ k+2} } \right)_{. {\kern 1pt} i} ({\rm {\bf
v}}^{(k+2)}_{.{\kern 1pt}  i} ) \in {\rm {\mathbb{H}}}^{n\times
n}$.  Taking into account Theorem
\ref{theor:char_polin} we obtain
\begin{equation}
\label{kyr19} \det \left( {\lambda{\rm {\bf I}} + {\rm {\bf
V}}^{ k+2} }
\right)_{.{\kern 1pt} i} \left( {{\rm {\bf
v}}^{(k+2)}_{.{\kern 1pt}
{\kern 1pt} i}}  \right) = d_{1} \lambda^{n - 1} + d_{2} \lambda^{n - 2} +
\ldots + d_{n},
\end{equation}
where $d_{s} = {\sum\limits_{\beta \in J_{s,\,n} {\left\{ {i}
\right\}}} {| \left( {{\rm {\bf V}}^{ k+2} } \right) _{\beta} ^{\beta} }| } $ is the sum of all principal
minors of order $s$ that contain the $i$-th column for all $s =
\overline {1,n - 1} $ and $d_{n} = \det \left( {{\rm {\bf V}}^{
k+2} } \right)$. Consequently we have ${\rm {\bf v}}^{(k+2)}_{.
{\kern 1pt} i} = \left( {{\begin{array}{*{20}c}
 {{\sum\limits_{l} {\bar{v}_{1l}^{(k)}  \bar{w}_{li}} } } \hfill \\
 {{\sum\limits_{l} {\bar{v}_{2l}^{(k)}  \bar{w}_{li}} } } \hfill \\
 { \vdots}  \hfill \\
 {{\sum\limits_{l} {\bar{v}_{nl}^{(k)}  \bar{w}_{li}} } } \hfill \\
\end{array}} } \right) = {\sum\limits_{l} {{\rm {\bf \bar{v}}}_{.\,l}^{(k)}  \bar{w}_{li}}
}$, where ${\rm {\bf \bar{v}}}_{.{\kern 1pt}  l}^{(k)}  $ is
the $l$th column-vector  of ${\rm {\bf \bar{V}}}^{k}=({\bf A}{\bf W})^{k}{\bf A}$  and ${\bf W}{\bf A}{\bf W}={\bf \bar{W}}=(\bar{w}_{li})$ for all
$l=\overline{1,n}$. Taking into account Theorem \ref{theorem:
determinant of hermitian matrix},  we obtain on the one hand
\begin{equation}
\label{kyr20}
\begin{array}{c}
  \det \left( {\lambda{\rm {\bf I}} + {\rm {\bf V}}^{ k+2} }
\right)_{.{\kern 1pt} i} \left( {\rm {\bf v}}^{(k+2)}_{.
{\kern 1pt} i} \right) = {\rm{cdet}} _{i} \left( {\lambda{\rm {\bf I}}
+ {\rm {\bf V}}^{ k+2} } \right)_{.{\kern 1pt} i}
\left( {\rm {\bf v}}^{(k+2)}_{.
{\kern 1pt} i}  \right) =\\
   = {\sum\limits_{l} {{\rm{cdet}} _{i} \left( {\lambda{\rm {\bf I}} + {\rm {\bf
   V}}^{k+2
}} \right)_{.{\kern 1pt} l} \left( {{\rm {\bf \bar{v}}}_{.\,l}^{(k)}  \bar{w}_{li}} \right) = {\sum\limits_{l}
{{\rm{cdet}} _{i} \left( {\lambda{\rm {\bf I}} + {\rm {\bf V}}^{ k+2} }
\right)_{.{\kern 1pt} i} \left( {{\rm {\bf \bar{v}}}_{.{\kern 1pt}
{\kern 1pt} l}^{(k)} }  \right) \cdot {\kern 1pt}} } }} \bar{w}_{li}
\end{array}
\end{equation}
On the other hand having changed the order of summation,  we get
for all $s = \overline {1,n - 1} $
\begin{equation}
\label{kyr21}
\begin{array}{c}
  d_{s} = {\sum\limits_{\beta \in J_{s,\,n} {\left\{ {i} \right\}}}
{\det \left( {{\rm {\bf V}}^{ k+2} } \right){\kern 1pt} {\kern
1pt} _{\beta} ^{\beta} } }  = {\sum\limits_{\beta \in J_{s,\,n}
{\left\{ {i} \right\}}} {{\rm{cdet}} _{i} \left( {{\rm {\bf V}}^{
k+2} } \right){\kern 1pt}
_{\beta} ^{\beta} } }  =\\
  {\sum\limits_{\beta \in J_{s,\,n}
{\left\{ {i} \right\}}} {{\sum\limits_{l} {{\rm{cdet}} _{i} \left(
{\left( {{\rm {\bf V}}^{ k+2} {\kern 1pt}} \right)_{.\, i} \left(
{{\rm {\bf \bar{v}}}_{.\, l}^{(k)} \bar{w}_{l\,i}} \right)} \right)}} }} \,
_{\beta} ^{\beta}=\\
  {\sum\limits_{l} { {{\sum\limits_{\beta \in J_{s,\,n} {\left\{ {i}
\right\}}} {{\rm{cdet}} _{i} \left( {\left( {{\rm {\bf V}}^{ k+2}
{\kern 1pt}}  \right)_{.{\kern 1pt} i} \left( {{\rm {\bf
\bar{v}}}_{.{\kern 1pt}  l}^{(k)} }  \right)} \right){\kern
1pt} _{\beta} ^{\beta} } }} }}  \cdot \bar{w}_{l{\kern 1pt} i}.
\end{array}
\end{equation}
By substituting (\ref{kyr20}) and (\ref{kyr21}) in (\ref{kyr19}),
and equating factors at $\bar{w} _ {l \, i} $ when $l = j $, we obtain
the equality (\ref{eq:char_cdet}).
$\blacksquare$

 By analogy can be proved
the following lemma.
\begin{lemma}
If ${\bf W}{\bf A}={\bf
U}=\left(u_{ij} \right)\in {\mathbb{H}}^{n\times n}$ is Hermitian with $ Ind{\kern 1pt} {\rm {\bf U}}=k$ and $t \in \mathbb{R}$, then
\[ {\rm{rdet}} _{j}  {( \lambda{\rm {\bf I}} + {\rm {\bf U}}^{k+2} )_{j\,.\,} ({\rm {\bf \bar{u}}}_{i.}^{(k)}
)}
   = r_{1}^{\left( {ij} \right)} \lambda^{n
- 1} +r_{2}^{\left( {ij} \right)} \lambda^{n - 2} + \ldots +
r_{n}^{\left( {ij} \right)},
\]
\noindent where  $r_{n}^{\left( {ij} \right)} = {\rm{rdet}} _{j}
{({\rm {\bf U}}^{k+2} )_{j\,.\,} ({\rm {\bf \bar{u}}}_{i.\,}^{(k)} )}$
and $r_{s}^{\left( {ij} \right)} = {{{\sum\limits_{\alpha \in
I_{s,n} {\left\{ {j} \right\}}} {{\rm{rdet}} _{j} \left( {({\rm
{\bf U}}^{k+2} )_{j\,.\,} ({\rm {\bf \bar{u}}}_{i.\,}^{(k)} )}
\right)\,_{\alpha} ^{\alpha} } }}}$ for all $s = \overline {1,n -
1} $ and $i,j = \overline {1,n}$.
\end{lemma}
\begin{theorem}\label{theor:det_rep_wdraz1}
If ${\rm {\bf A}} \in  {\mathbb{H}}^{m\times n}$, ${\bf W}\in {\mathbb{H}}^{n\times m}$, and ${\bf A}{\bf W}={\bf
V}=\left(v_{ij} \right)\in {\mathbb{H}}^{m\times m}$ is Hermitian with $k= {\rm max}\{Ind({\bf A}{\bf W}), Ind({\bf W}{\bf A})\}$ and $\rank({\bf A}{\bf W})^{k+1} =
\rank({\bf A}{\bf W})^{k} = r$, then the W-weighted Drazin inverse ${\rm {\bf
A}}_{d,W} = \left( {a_{ij}^{d,W} } \right) \in {\rm
{\mathbb{H}}}^{m\times n} $ with respect to ${\bf W}$ possess the following determinantal
representations:
\begin{equation}
\label{eq:dr_rep_wcdet} a_{ij}^{d,W}  = {\frac{{{\sum\limits_{\beta
\in J_{r,\,m} {\left\{ {i} \right\}}} {{\rm{cdet}} _{i} \left(
{\left( {{\rm {\bf A}}{\bf W}} \right)^{k+2}_{\,. \,i} \left( {{\rm {\bf
\bar{v}}}_{.j}^{ (k)} }  \right)} \right){\kern 1pt} {\kern 1pt} _{\beta}
^{\beta} } } }}{{{\sum\limits_{\beta \in J_{r,\,\,m}} {{\left|
{\left( { {\bf A}{\bf W}} \right)^{k+2}{\kern 1pt} _{\beta} ^{\beta}
}  \right|}}} }}},
\end{equation}
where ${\rm {\bf \bar{v}}}_{.j}^{(k)} $ is the $j$th column of  ${\rm
{\bf \bar{V}}}^{k}=({\bf A}{\bf W})^{k}{\bf A} $ for all $j=\overline{1,m}$.
\end{theorem}
{\bf Proof}.   The
matrix $\left( {\lambda {\rm {\bf I}}_m + ({\rm {\bf A}}{\bf W})^{k + 2}} \right)^{
- 1} \in {\rm {\mathbb{H}}}^{m\times m}$ is a full-rank
Hermitian matrix. Taking into account Theorem \ref{inver_her} it
has an inverse, which we represent as a left inverse matrix
\[
\left( {\lambda {\rm {\bf I}}_m + ({\rm {\bf A}}{\bf W})^{k + 2}} \right)^{
- 1}
= {\frac{{1}}{{\det \left( {\lambda {\rm {\bf I}}_m + ({\rm {\bf A}}{\bf W})^{k + 2}} \right)}}}\left( {{\begin{array}{*{20}c}
 {L_{11}}  \hfill & {L_{21}}  \hfill & {\ldots}  \hfill & {L_{m1}}  \hfill
\\
 {L_{12}}  \hfill & {L_{22}}  \hfill & {\ldots}  \hfill & {L_{m2}}  \hfill
\\
 {\ldots}  \hfill & {\ldots}  \hfill & {\ldots}  \hfill & {\ldots}  \hfill
\\
 {L_{1m}}  \hfill & {L_{2m}}  \hfill & {\ldots}  \hfill & {L_{mm}}  \hfill
\\
\end{array}} } \right),
\]
\noindent where $L_{ij} $ is a left $ij$-th cofactor  of a matrix
${\lambda {\rm {\bf I}}_m + ({\rm {\bf A}}{\bf W})^{k + 2}}$. Then we have
\[\begin{array}{l}
   \left( {\lambda {\rm {\bf I}}_m + ({\rm {\bf A}}{\bf W})^{k + 2}} \right)^{
- 1}({\rm {\bf A}}{\bf W})^{k}{\bf A}  = \\
  ={\frac{{1}}{{\det \left( {\lambda {\rm {\bf I}}_m + ({\rm {\bf A}}{\bf W})^{k + 2}} \right)}}}\left(
{{\begin{array}{*{20}c}
 {{\sum\limits_{s = 1}^{m} {L_{s1} \bar{v}_{s1}^{ (k)} } } } \hfill &
{{\sum\limits_{s = 1}^{m} {L_{s1} \bar{v}_{s2}^{ (k)} } } } \hfill &
{\ldots}
\hfill & {{\sum\limits_{s = 1}^{m} {L_{s1} \bar{v}_{sn}^{ (k)} } } } \hfill \\
 {{\sum\limits_{s = 1}^{m} {L_{s2} \bar{v}_{s1}^{ (k)} } } } \hfill &
{{\sum\limits_{s = 1}^{m} {L_{s2} \bar{v}_{s2}^{ (k)} } } } \hfill &
{\ldots}
\hfill & {{\sum\limits_{s = 1}^{m} {L_{s2} \bar{v}_{sn}^{ (k)} } } } \hfill \\
 {\ldots}  \hfill & {\ldots}  \hfill & {\ldots}  \hfill & {\ldots}  \hfill
\\
 {{\sum\limits_{s = 1}^{m} {L_{sm} \bar{v}_{s1}^{ (k)} } } } \hfill &
{{\sum\limits_{s = 1}^{m} {L_{sm} \bar{v}_{s2}^{ (k)} } } } \hfill &
{\ldots}
\hfill & {{\sum\limits_{s = 1}^{m} {L_{sm} \bar{v}_{sn}^{ (k)} } } } \hfill \\
\end{array}} } \right).
\end{array}
\]
By
(\ref{eq:lim_repr_AW}) and using the definition of a left cofactor, we obtain
\begin{equation}
\label{eq:A_d1}  {\rm {\bf A}}_{ d,W}   = {\mathop {\lim}
\limits_{\alpha \to 0}} \left( {{\begin{array}{*{20}c}
 {{\frac{{{\rm cdet} _{1} \left( {\lambda {\rm {\bf I}}_m + ({\rm {\bf A}}{\bf W})^{k + 2}} \right)_{.1} \left( {{\rm {\bf \bar{v}}}_{.1}^{ (k)} }
\right)}}{{\det \left( {\lambda {\rm {\bf I}}_m + ({\rm {\bf A}}{\bf W})^{k + 2}} \right)}}}} \hfill & {\ldots}  \hfill & {{\frac{{{\rm
cdet} _{1} \left( {\lambda {\rm {\bf I}}_m + ({\rm {\bf A}}{\bf W})^{k + 2}} \right)_{.1} \left( {{\rm {\bf \bar{v}}}_{.n}^{(k)} } \right)}}{{\det
\left( {\lambda {\rm {\bf I}}_m + ({\rm {\bf A}}{\bf W})^{k + 2}} \right)}}}} \hfill \\
 {\ldots}  \hfill & {\ldots}  \hfill & {\ldots}  \hfill \\
 {{\frac{{{\rm cdet} _{n} \left( {\lambda {\rm {\bf I}}_m + ({\rm {\bf A}}{\bf W})^{k + 2}} \right)_{.n} \left( {{\rm {\bf \bar{v}}}_{.1}^{ (k)} }
\right)}}{{\det \left( {\lambda {\rm {\bf I}}_m + ({\rm {\bf A}}{\bf W})^{k + 2}} \right)}}}} \hfill & {\ldots}  \hfill & {{\frac{{{\rm
cdet} _{n} \left( {\lambda {\rm {\bf I}}_m + ({\rm {\bf A}}{\bf W})^{k + 2}} \right)_{.m} \left( {{\rm {\bf \bar{v}}}_{.n}^{(k)} } \right)}}{{\det
\left( {\lambda {\rm {\bf I}}_m + ({\rm {\bf A}}{\bf W})^{k + 2}} \right)}}}} \hfill \\
\end{array}} } \right).
\end{equation}
By Theorem \ref{theor:char_polin} we have
 \[\det \left( {\lambda {\rm {\bf I}}_m + ({\rm {\bf A}}{\bf W})^{k + 2}} \right) = \lambda ^{m} + d_{1}
\lambda ^{m - 1} + d_{2} \lambda ^{m - 2} + \ldots + d_{m},\]
where
$d_{s} = {\sum\limits_{\beta \in J_{s,\,m}}  {{\left| {\left( {\lambda {\rm {\bf I}}_m + ({\rm {\bf A}}{\bf W})^{k + 2}} \right){\kern 1pt} {\kern 1pt} _{\beta}
^{\beta} } \right|}}} $ is a sum of principal minors of $({\rm {\bf A}}{\bf W})^{k + 2}$ of order $s$ for all  $s = \overline {1,m - 1} $ and
$d_{m} = \det( {\rm {\bf A}}{\bf W})^{k + 2}$.
Since $\rank({\rm {\bf A}}{\bf W})^{k + 2}=\rank({\rm {\bf A}}{\bf W})^{k+1} = \rank({\rm {\bf A}}{\bf W})^{k} = r$,
then $d_{m} = d_{m - 1} = \ldots = d_{r + 1} = 0$. It follows that
$\det \left( {\lambda {\rm {\bf I}}_m + ({\rm {\bf A}}{\bf W})^{k + 2}} \right)
= \lambda ^{m} + d_{1} \lambda ^{m - 1} + d_{2}\lambda ^{m - 2} +
\ldots + d_{r} \lambda ^{m - r}$.
Using (\ref{eq:char_cdet})  we have
 \[{\rm{cdet}} _{i} \left(\lambda {\rm {\bf I}}_m + {({\rm {\bf A}}{\bf W})^{ k+2} }
\right)_{.\,i} \left( { \bar{{\bf v}}_{.j}^{ (k)} }  \right) = c_{1}^{\left( {ij} \right)} \lambda^{m - 1}
+ c_{2}^{\left( {ij} \right)} \lambda^{m - 2} + \ldots + c_{m}^{\left(
{ij} \right)} \]
for  $i = \overline {1,m} $ and $j = \overline {1,n} $, where $c_{s}^{\left( {ij}
\right)} = {\sum\limits_{\beta \in J_{s,\,m} {\left\{ {i}
\right\}}} {{\rm{cdet}} _{i} \left( {({\rm {\bf A}}{\bf W})^{ k+2} _{.\,i}
\left( {{\rm {\bf \bar{v}}}_{.j}^{ (k)} } \right)} \right){\kern 1pt}
{\kern 1pt} _{\beta }^{\beta} } } $ for all $s = \overline {1,m -
1} $ and $c_{m}^{\left( {ij} \right)} = {\rm{cdet}} _{i} ({\rm {\bf A}}{\bf W})^{ k+2} _{.i} \left( {{\rm {\bf \bar{v}}}_{.j}^{ (k)}
} \right)$.
We shall prove that $c_{k}^{\left( {ij} \right)} = 0$, when $k \ge r +
1$ for   $i = \overline {1,m} $ and $j = \overline {1,n} $.
By Lemma \ref{lem:rank_col}
$\left( {({\rm {\bf A}}{\bf W})^{ k+2} _{.\,i}
\left( {{\rm {\bf \bar{v}}}_{.j}^{ (k)} } \right)} \right) \le r$, then the matrix
$\left( {({\rm {\bf A}}{\bf W})^{ k+2} _{.\,i}
\left( {{\rm {\bf \bar{v}}}_{.j}^{ (k)} } \right)} \right)$ has no more $r$ right-linearly
independent columns.
Consider $\left( {({\rm {\bf A}}{\bf W})^{ k+2} _{.\,i}
\left( {{\rm {\bf \bar{v}}}_{.j}^{ (k)} } \right)} \right) {\kern 1pt} _{\beta}
^{\beta}  $, when $\beta \in J_{s,m} {\left\{ {i} \right\}}$. It
is a principal submatrix of  $\left( {({\rm {\bf A}}{\bf W})^{ k+2} _{.\,i}
\left( {{\rm {\bf \bar{v}}}_{.j}^{ (k)} } \right)} \right)$ of order
$s \ge r + 1$. Deleting both its $i$-th row and column, we obtain
a principal submatrix of order $s - 1$ of  $({\rm {\bf A}}{\bf W})^{ k+2}$.
 We denote it by ${\rm {\bf M}}$. The following cases are
possible.
\begin{itemize}
  \item Let $s = r + 1$ and $\det {\rm {\bf M}} \ne 0$. In this case all
columns of $ {\rm {\bf M}} $ are right-linearly independent. The
addition of all of them on one coordinate to columns of
 $\left( {({\rm {\bf A}}{\bf W})^{ k+2} _{.\,i}
\left( {{\rm {\bf \bar{v}}}_{.j}^{ (k)} } \right)} \right)
{\kern 1pt} _{\beta} ^{\beta}$ keeps their right-linear
independence. Hence, they are basis in a matrix $\left( {({\rm {\bf A}}{\bf W})^{ k+2} _{.\,i}
\left( {{\rm {\bf \bar{v}}}_{.j}^{ (k)} } \right)} \right) {\kern 1pt} _{\beta}
^{\beta} $, and  the $i$-th
column is the right linear combination of its basis columns. From
this by Theorem \ref{theorem:colum_combin}, we get ${\rm{cdet}}
_{i} \left( {({\rm {\bf A}}{\bf W})^{ k+2} _{.\,i}
\left( {{\rm {\bf \bar{v}}}_{.j}^{ (k)} } \right)} \right) {\kern 1pt}
_{\beta} ^{\beta}  = 0$, when $\beta \in J_{s,n} {\left\{ {i}
\right\}}$ and $s= r + 1$.
  \item If $s = r + 1$ and $\det {\rm {\bf M}} = 0$, than $p$, ($p < s$),
columns are basis in  ${\rm {\bf M}}$ and in  $\left( {({\rm {\bf A}}{\bf W})^{ k+2} _{.\,i}
\left( {{\rm {\bf \bar{v}}}_{.j}^{ (k)} } \right)} \right) {\kern 1pt} _{\beta} ^{\beta} $.
Then by Theorem \ref{theorem:colum_combin}, ${\rm{cdet}} _{i} \left( {({\rm {\bf A}}{\bf W})^{ k+2} _{.\,i}
\left( {{\rm {\bf \bar{v}}}_{.j}^{ (k)} } \right)} \right){\kern 1pt} _{\beta}
^{\beta}  = 0$ as well.
  \item If $s > r + 1$, then  $\det {\rm {\bf M}} = 0$ and
$p$, ($p < r$), columns are basis in the both matrices ${\rm {\bf
M}}$ and $\left( {({\rm {\bf A}}{\bf W})^{ k+2} _{.\,i}
\left( {{\rm {\bf \bar{v}}}_{.j}^{ (k)} } \right)} \right)
{\kern 1pt} _{\beta} ^{\beta}  $. Then by Theorem
\ref{theorem:colum_combin}, we also
have ${\rm{cdet}} _{i} \left( {({\rm {\bf A}}{\bf W})^{ k+2} _{.\,i}
\left( {{\rm {\bf \bar{v}}}_{.j}^{ (k)} } \right)} \right) {\kern 1pt} _{\beta} ^{\beta}  = 0.$
\end{itemize}
Thus in all cases we have ${\rm{cdet}} _{i} \left( {({\rm {\bf A}}{\bf W})^{ k+2} _{.\,i}
\left( {{\rm {\bf \bar{v}}}_{.j}^{ (k)} } \right)} \right) {\kern 1pt} _{\beta} ^{\beta}  = 0$,
when $\beta \in J_{s,m} {\left\{ {i} \right\}}$ and $r + 1 \le s <
m$. From here if $r + 1 \le s < m$, then
\[c_{s}^{\left( {ij} \right)} = {\sum\limits_{\beta \in J_{s,\,m}
{\left\{ {i} \right\}}} {{\rm{cdet}} _{i} \left( {({\rm {\bf A}}{\bf W})^{ k+2} _{.\,i}
\left( {{\rm {\bf \bar{v}}}_{.j}^{ (k)} } \right)} \right) {\kern 1pt} _{\beta} ^{\beta} } }=0,\] and
$c_{m}^{\left( {ij} \right)} = {\rm{cdet}} _{i} \left( {({\rm {\bf A}}{\bf W})^{ k+2} _{.\,i}
\left( {{\rm {\bf \bar{v}}}_{.j}^{ (k)} } \right)} \right) = 0$ for all $i,j = \overline {1,n} $.

Hence, ${\rm{cdet}} _{i} \left( {\lambda {\rm {\bf I}} + ({\rm {\bf A}}{\bf W})^{ k+2}} \right)_{.\,i} \left( {{\rm {\bf \bar{v}}}_{.\,j}^{ (k)} }
\right) =c_{1}^{\left( {ij} \right)} \lambda ^{m - 1} +
c_{2}^{\left( {ij} \right)} \lambda ^{m - 2} + \ldots +
c_{r}^{\left( {ij} \right)} \lambda ^{m - r}$ for  $i = \overline {1,m} $ and $j = \overline {1,n} $. By substituting these values in the matrix from
(\ref{eq:A_d1}), we obtain
\[\begin{array}{c}
  {\rm {\bf A}}_{ d,W}  = {\mathop {\lim} \limits_{\lambda \to 0}} \left(
{{\begin{array}{*{20}c}
 {{\frac{{c_{1}^{\left( {11} \right)} \lambda ^{m - 1} + \ldots +
c_{r}^{\left( {11} \right)} \lambda ^{m - r}}}{{\lambda ^{m} + d_{1}
\lambda ^{m - 1} + \ldots + d_{r} \lambda ^{m - r}}}}} \hfill &
{\ldots}  \hfill & {{\frac{{c_{1}^{\left( {1n} \right)} \lambda ^{m
- 1} + \ldots + c_{r}^{\left( {1n} \right)} \lambda ^{m -
r}}}{{\lambda ^{m} + d_{1} \lambda
^{m - 1} + \ldots + d_{r} \lambda ^{m - r}}}}} \hfill \\
 {\ldots}  \hfill & {\ldots}  \hfill & {\ldots}  \hfill \\
 {{\frac{{c_{1}^{\left( {m1} \right)} \lambda ^{m - 1} + \ldots +
c_{r}^{\left( {m1} \right)} \lambda ^{m - r}}}{{\lambda ^{m} + d_{1}
\lambda ^{m - 1} + \ldots + d_{r} \lambda ^{m - r}}}}} \hfill &
{\ldots}  \hfill & {{\frac{{c_{1}^{\left( {mn} \right)} \lambda ^{m
- 1} + \ldots + c_{r}^{\left( {mn} \right)} \lambda ^{m -
r}}}{{\lambda ^{m} + d_{1} \lambda
^{m - 1} + \ldots + d_{r} \lambda ^{m - r}}}}} \hfill \\
\end{array}} } \right) =\\
  \left( {{\begin{array}{*{20}c}
 {{\frac{{c_{r}^{\left( {11} \right)}} }{{d_{r}} }}} \hfill & {\ldots}
\hfill & {{\frac{{c_{r}^{\left( {1n} \right)}} }{{d_{r}} }}} \hfill \\
 {\ldots}  \hfill & {\ldots}  \hfill & {\ldots}  \hfill \\
 {{\frac{{c_{r}^{\left( {m1} \right)}} }{{d_{r}} }}} \hfill & {\ldots}
\hfill & {{\frac{{c_{r}^{\left( {mn} \right)}} }{{d_{r}} }}} \hfill \\
\end{array}} } \right).
\end{array}
\]
Here $c_{r}^{\left( {ij} \right)} = {\sum\limits_{\beta \in
J_{r,\,m} {\left\{ {i} \right\}}} {{\rm{cdet}} _{i} \left( {\left(
{{\rm {\bf A}}^{k+1}} \right)_{\,.\,i} \left( {{\rm {\bf
a}}_{.j}^{ (k)} }  \right)} \right){\kern 1pt} {\kern 1pt} _{\beta}
^{\beta} } } $ and $d_{r} = {\sum\limits_{\beta \in J_{r,\,m}}
{{\left| {\left(
{{\rm {\bf A}}^{k+1}} \right){\kern
1pt} {\kern 1pt} _{\beta} ^{\beta} } \right|}}} $. Thus, we have
obtained the determinantal representation of ${\rm {\bf A}}_{d,W}  $ by (\ref{eq:dr_rep_wcdet}).
$\blacksquare$

By analogy can be proved the following theorem.
\begin{theorem}\label{theor:det_rep_wdraz2}
If ${\rm {\bf A}} \in  {\mathbb{H}}^{m\times n}$, ${\bf W}\in {\mathbb{H}}^{n\times m}$, and ${\bf W}{\bf A}={\bf
U}=\left(u_{ij} \right)\in {\mathbb{H}}^{n\times n}$ is Hermitian with $k= {\rm max}\{Ind({\bf A}{\bf W}), Ind({\bf W}{\bf A})\}$ and $\rank({\bf W}{\bf A})^{k+1} =
\rank({\bf W}{\bf A})^{k} = r$, then the W-weighted Drazin inverse ${\rm {\bf
A}}_{d,W} = \left( {a_{ij}^{d,W} } \right) \in {\rm
{\mathbb{H}}}^{m\times n} $ with respect to ${\bf W}$ possess the following determinantal
representations:
\begin{equation}
\label{eq:dr_rep_wrdet} a_{ij}^{d,W}  = {\frac{{{\sum\limits_{\alpha
\in I_{r,n} {\left\{ {j} \right\}}} {{\rm{rdet}} _{j} \left(
{({\bf W}{\rm {\bf A}} )^{ k+2}_{j\,.\,} ({\rm {\bf \bar{u}}}_{i.\,}^{ (k)} )}
\right)\,_{\alpha} ^{\alpha} } }}}{{{\sum\limits_{\alpha \in
I_{r,\,n}}  {{\left| {\left({\bf W} {{\rm {\bf A}} } \right)^{k+2}{\kern
1pt}  _{\alpha} ^{\alpha} } \right|}}} }}}.
\end{equation}
where  ${\rm {\bf
\bar{u}}}_{i.}^{(k)} $ is  the $i$th row of  ${\rm
{\bf \bar{U}}}^{k}={\bf A}({\bf W}{\bf A})^{k} $ for all $i=\overline{1,
n}$.
\end{theorem}
\section{An example}
In this section, we give an example to illustrate our results. Let
us consider the matrices
\[{\bf A}=\begin{pmatrix}
  0 & i & 0 \\
  k & 1 & i \\
 1 & 0 & 0\\
  1 & -k & -j
\end{pmatrix},\,\, {\bf W}=\begin{pmatrix}
 k & 0 & i & 0 \\
 -j & k & 0 & 1 \\
 0 & 1 & 0 & -k
\end{pmatrix}.\]
Then
\[{\bf V}={\bf A}{\bf W}=\begin{pmatrix}
 -k & -j & 0 & i \\
 -1-j & i+k & j & 1+j \\
 k & 0 & i & 0\\
 -i+k & 1-j & i & i-k
\end{pmatrix},\,\,{\bf U}={\bf W}{\bf A}=\begin{pmatrix}
  i & j & 0 \\
  0 & k & 0 \\
 0 & 0 & 0
\end{pmatrix},\]
and $\rank{\bf W}=3$, $\rank{\bf V}=3$, $\rank{\bf V}^{3}=\rank {\bf V}^{2}=2$, $\rank{\bf U}^{2}=\rank {\bf U}=2$.
Therefore, $ Ind\,{\bf V}=2$, $ Ind\,{\bf U}=1$, and $k= {\rm max}\{Ind({\bf A}{\bf W}), Ind({\bf W}{\bf A})\}=2$.
It's evident that obtaining the W-weighted Drazin inverse of ${\bf A}$ with respect to ${\bf W}$ by using the matrix ${\bf U}$ by (\ref{eq:det_repr_u wdraz}) is more convenient.
We have \[{\bf U}^{2}=\begin{pmatrix}
  -1 & i+k & 0 \\
  0 & -1 &  \\
 0 & 0 & 0
\end{pmatrix},\,\,\,{\bf U}^{5}=\begin{pmatrix}
  i & 2+3j & 0 \\
  0 & k &  \\
 0 & 0 & 0
\end{pmatrix},\]\[({\bf U}^{5})^{*}=\begin{pmatrix}
  -i & 0 & 0 \\
  2-3j & -k &  \\
 0 & 0 & 0
\end{pmatrix},\left({\bf U}^{5}\right)^{*}{\bf U}^{5}=\begin{pmatrix}
  1 & -2i-3k & 0 \\
  2i+3k & 14 &  \\
 0 & 0 & 0
\end{pmatrix},\] \[
\hat{{\rm {\bf U}}}=({\bf U}^{ 5})^{*}{\bf U}^{2}=\begin{pmatrix}
  i & 1+j & 0 \\
  -2+3j & -i+6k &  \\
 0 & 0 & 0
\end{pmatrix},{\bf W}^{*}=\begin{pmatrix}
 -k & j & 0  \\
 0 & -k & 1 \\
 -i & 0 &  0\\
 0 & 1 & k
\end{pmatrix},\]\[
{\bf W}^{*}{\bf W}=\begin{pmatrix}
 2 & i & -j & j  \\
 -i & 2 & 0 & -2k \\
 j & 0 &  1 & 0\\
  -j & 2k & 0 & 2 \\
\end{pmatrix},
\hat{{\rm {\bf W}}}={\bf W}^{*}{\bf U}^{2}=\begin{pmatrix}
 -k & 1-2j & 0  \\
 0 & i+k & 0 \\
 i & 1+j &  0\\
 0 & -1 & 0
\end{pmatrix}.
\]
Since by (\ref{eq:det_repr_u wdraz})
\begin{multline*}
 {a}_{11}^{d,{\bf W}}=\\
\frac{\sum\limits_{t = 1}^{3}
{{\sum\limits_{\beta \in I_{3,\,4} {\left\{ {1}
\right\}}} {{\rm{cdet}} _{1} \left( {\left(   {\bf W}^{*} {\bf W}\right)_{.1 } (\hat{{\rm {\bf w}}}_{.t})} \right) {\kern 1pt} _{\beta} ^{\beta} } }
}
{{\sum\limits_{\beta \in J_{2,\,3} {\left\{ {t}
\right\}}} {{\rm{cdet}} _{t} \left( {\left( { \left({\bf U}^{ 5} \right)^{*}{\bf U}^{ 5}
} \right)_{.t } (\hat{{\rm {\bf u}}}_{.1})} \right) {\kern 1pt} _{\beta} ^{\beta} } }
}}
{{\sum\limits_{\beta \in J_{3,\,4}} {{\left| {\left(  {\bf W}^{*}{\bf W}\right){\kern 1pt} _{\beta} ^{\beta}
}  \right|}}}{\sum\limits_{\beta \in J_{2,\,
3}} {{\left| {\left( { \left({\bf U}^{ 5} \right)^{*}{\bf U}^{ 5}
} \right){\kern 1pt} _{\beta} ^{\beta}
}  \right|}}}},\end{multline*}
where
\begin{multline*}{{\sum\limits_{\beta \in I_{3,\,4} {\left\{ {1}
\right\}}} {{\rm{cdet}} _{1} \left( {\left(   {\bf W}^{*} {\bf W}\right)_{.1 } (\hat{{\rm {\bf w}}}_{.1})} \right) {\kern 1pt} _{\beta} ^{\beta} } }
}=\\{\rm{cdet}} _{1}\begin{pmatrix}
 k & i & -j  \\
 0 & 2 & 0 \\
 i & 0 &  1
\end{pmatrix}+{\rm{cdet}} _{1}\begin{pmatrix}
 k & i & j  \\
 0 & 2 & -2k \\
 0 & 2k &  1
\end{pmatrix}+{\rm{cdet}} _{1}\begin{pmatrix}
 k & -j & j  \\
 i & 1 & 0 \\
 0 & 0 &  2
\end{pmatrix}=0, \end{multline*}
\begin{multline*}{{\sum\limits_{\beta \in I_{3,\,4} {\left\{ {1}
\right\}}} {{\rm{cdet}} _{1} \left( {\left(   {\bf W}^{*} {\bf W}\right)_{.1 } (\hat{{\rm {\bf w}}}_{.2})} \right) {\kern 1pt} _{\beta} ^{\beta} } }
}=-2j,
{{\sum\limits_{\beta \in I_{3,\,4} {\left\{ {1}
\right\}}} {{\rm{cdet}} _{1} \left( {\left(   {\bf W}^{*} {\bf W}\right)_{.1 } (\hat{{\rm {\bf w}}}_{.3})} \right) {\kern 1pt} _{\beta} ^{\beta} } }}=0,\\
{\sum\limits_{\beta \in J_{3,\,4}} {{\left| {\left(  {\bf W}^{*}{\bf W}\right){\kern 1pt} _{\beta} ^{\beta}
}  \right|}}}=2,\end{multline*}
and
\begin{multline*}\sum\limits_{\beta \in J_{2,\,3} {\left\{ {1}
\right\}}}{{\rm{cdet}} _{1} \left( {\left( { \left({\bf U}^{ 5} \right)^{*}{\bf U}^{ 5}
} \right)_{.1 } (\hat{{\rm {\bf u}}}_{.1})} \right) {\kern 1pt} _{\beta} ^{\beta} }=\\{\rm{cdet}} _{1}\begin{pmatrix}
 i & -2i-3k   \\
 -2+3j & 14
\end{pmatrix}+{\rm{cdet}} _{1}\begin{pmatrix}
 i & 0   \\
 0 & 0
\end{pmatrix}=i,\end{multline*}
\begin{multline*}\sum\limits_{\beta \in J_{2,\,3} {\left\{ {2}
\right\}}}{{\rm{cdet}} _{2} \left( {\left( { \left({\bf U}^{ 5} \right)^{*}{\bf U}^{ 5}
} \right)_{.2 } (\hat{{\rm {\bf u}}}_{.1})} \right) {\kern 1pt} _{\beta} ^{\beta} }=0,\\\sum\limits_{\beta \in J_{2,\,3} {\left\{ {3}
\right\}}}{{\rm{cdet}} _{3} \left( {\left( { \left({\bf U}^{ 5} \right)^{*}{\bf U}^{ 5}
} \right)_{.3 } (\hat{{\rm {\bf u}}}_{.1})} \right) {\kern 1pt} _{\beta} ^{\beta} }=0,{\sum\limits_{\beta \in J_{2,\,
3}} {{\left| {\left( { \left({\bf U}^{ 5} \right)^{*}{\bf U}^{ 5}
} \right){\kern 1pt} _{\beta} ^{\beta}
}  \right|}}}=1,\end{multline*}
then
\[
 {a}_{11}^{d,{\bf W}}=
\frac{(0\cdot i)+(-2j\cdot 0)+(0\cdot 0)}{2\cdot 1 }=0.
\]
Continuing in the same way, we finally get,
 \begin{equation}
\label{ex:wdr_rep1}{\rm {\bf A}}_{d,W}= \begin{pmatrix}
  0 & -i & 0 \\
  0 &0 &0  \\
 -1 & 5i-2k & 0\\
 0 & 0 & 0
\end{pmatrix}.\end{equation}
By (\ref{eq:det_repr_A*A}) we can obtain
\[\left({\bf U}^{5}\right)^{+}=\begin{pmatrix}
  -i & -3+2j & 0 \\
  0 & -k &  \\
 0 & 0 & 0
\end{pmatrix},\,\,({\bf W}{\bf A})^{D}={\bf U}^{2}\left({\bf U}^{5}\right)^{+}{\bf U}^{2}=\begin{pmatrix}
  -i & -5 & 0 \\
  0 & -k &  \\
 0 & 0 & 0
\end{pmatrix}.\]
We verify (\ref{ex:wdr_rep1}) by (\ref{eq:ADW}). Indeed,
\[{\bf W}{\rm {\bf A}}_{d,{\bf W}}=\begin{pmatrix}
 k & 0 & i & 0 \\
 -j & k & 0 & 1 \\
 0 & 1 & 0 & -k
\end{pmatrix}\begin{pmatrix}
  0 & -i & 0 \\
  0 &0 &0  \\
 -1 & 5i-2k & 0\\
 0 & 0 & 0
\end{pmatrix}=\begin{pmatrix}
  -i & -5 & 0 \\
  0 & -k &  \\
 0 & 0 & 0
\end{pmatrix}=({\bf W}{\bf A})^{D}.\]
We also obtain  the W-weighted Drazin inverse of ${\bf A}$ with respect to ${\bf W}$ by (\ref{eq:det_rep1_wdraz}), then we have
   \begin{equation}\label{ex:wdr_rep2}
{\rm {\bf A}}_{d,{\bf W}}={\bf A}\left(({\bf W}{\bf A})^{D}) \right)^{2}=\begin{pmatrix}
  0 & -i & 0 \\
  -k &6+5i &0  \\
 -1 & 5i+5k & 0\\
 -1 & 5i+6k & 0
\end{pmatrix},
\end{equation}
The W-weighted Drazin inverse in (\ref{ex:wdr_rep2}) different from (\ref{ex:wdr_rep1}). It can be explained  that the Jordan normal form of ${\bf W}{\bf A}$ is unique only up to the order of the Jordan blocks.
We get their complete equality, if ${\rm {\bf A}}_{d,{\bf W}}$ from (\ref{ex:wdr_rep2}) be  left-multiply by the nonsingular matrix ${\bf P}$ which is  the product of multiplication of the following elementary matrices,
\[{\bf P}={\bf P}_{2,4}(-k)\cdot{\bf P}_{4,3}(-1)\cdot{\bf P}_{3,4}(-6)\cdot {\bf P}_{4,1}(-j)=\begin{pmatrix}
 1 & 0 & 0 & 0 \\
  0 &1 &0 &-k  \\
 0 & 0 & 7 & -6\\
 -j &0 & -1 & 1
\end{pmatrix}.\]
Note that we used  Maple with the package CLIFFORD in the calculations.


\begin{thebibliography}{40}
\bibitem{cl} R.E. Cline, T.N.E. Greville, A Drazin inverse for rectangular matrices., Linear Algebra
Appl. 29 (1980) 53--62.
\bibitem{wei1} Y. Wei, Integral representation of the W-weighted Drazin inverse, Appl. Math. Comput. 144 (2003) 3--10.
\bibitem{wei2} Y. Wei, C.-W. Woo, T. Lei, A note on the perturbation of the W-
weighted Drazin Inverse, Appl. Math. Comput.149 (2004)  423--430.
\bibitem {wei4}Y. Wei, A characterization for the W -weighted Drazin inverse and a Cramer rule for the W -
weighted Drazin inverse solution, Appl. Math. Comput. 125 (2002) 303--310.
\bibitem{ze} Z. Al-Zhour, A. Kili\c{c}man, M. H. Abu Hassa, New representations for weighted
Drazin inverse of matrices, Int. Journal of Math. Analysis 1(15) (2007) 697--708.
\bibitem{mos} D. Mosi\'{c}, D. S. Djordjevi\'{c}, Additive results for theWg-Drazin inverse, Linear
Algebra Appl. 432 (2010) 2847-2860.
 \bibitem{ky1} I. Kyrchei, Cramer's rule for quaternion
 systems of linear equations, Journal of Mathematical Sciences 155 (6) (2008) 839-858.
 \bibitem{ky2}   I. Kyrchei, The theory of the column and row determinants in a quaternion linear algebra,  in: Albert R. Baswell (Eds.), Advances in Mathematics Research 15,  Nova Sci. Publ., New York, 2012, pp. 301-359.
     \bibitem{di} J. Dieudonne, Les determinants sur un corps non-commutatif, Bull. Soc.
Math. France 71 (1943), 27--45.
     \bibitem{stu} E. Study, Zur Theorie der linearen Gleichungen, Acta Math. 42 (1920)
1--61.
    \bibitem{mo} E. H. Moore, On the determinant of an hermitian matrix of
quaternionic elements, Bull. Amer. Math. Soc. 28 (1922) 161--162.
 \bibitem{dy} F. J. Dyson, Quaternion determinants, Helv. Phys. Acta 45 (1972), 289--302.
   \bibitem{ch} L. Chen, Definition of determinant and Cramer solution over the
quaternion field, Acta Math. Sinica (N.S.) 7 (1991) 171--180.
     \bibitem{ge}I. M. Gelfand, V. S. Retakh, Determinants of matrices over noncommutative
rings, Functional Anal. Appl. 25 (1991) 91--102.
\bibitem{ky3} I. Kyrchei, Determinantal representations of the Moore-Penrose inverse over the quaternion skew field and corresponding Cramer's rules, Linear Multilinear Algebra 59 (2011) 413-431.
     \bibitem{ky33} I. Kyrchei, Determinantal representation of the Moore–Penrose inverse matrix over the quaternion skew field, Journal of Mathematical Sciences 180(1) (2012)  23--33.
   \bibitem{ky4} I. Kyrchei, Determinantal representations of the Drazin inverse over the quaternion skew field with applications to some matrix equations, Appl. Math. Comput. 238 (2014), pp. 193-207.

  \bibitem{song1}G. Song, Determinantal representations of the generalized inverses $A_ {T, S}^{(2)}$ over the quaternion skew field with applications, Journal of Applied Mathematics and Computing 39 (2012) 201--220.
        \bibitem{song3}G. Song, Q. Wang, H. Chang. Cramer rule for the unique solution of restricted matrix
equations over the quaternion skew field. Comput. Math. Appl. 61 (2011)1576--1589.
        \bibitem{song2}G. Song,
     Characterization of the W-weighted Drazin inverse over the quaternion skew field with applications, Electronic Journal of Linear Algebra 26 (2013) 1--14.
\bibitem{ca1} S. L. Campbell and C.D. Meyer, Generalized inverse of linear transformations, Corrected reprint of the 1979 original. Dover Publications, Inc., New York, 1991.
\bibitem{hu}
L. Huang, W. So, On left eigenvalues of a quaternionic matrix, Linear
Algebra Appl. 323 (2001) 105-116.
 \bibitem{so}
W. So,  Quaternionic left eigenvalue problem, Southeast Asian Bulletin of
Mathematics 29 (2005)  555-565.
 \bibitem{wo}
R. M. W. Wood, Quaternionic eigenvalues, Bull. Lond. Math. Soc. 17
(1985)137-138.
\bibitem{br}
J.L. Brenner, Matrices of quaternions, Pac. J. Math. 1 (1951) 329-335.
 \bibitem{ma} E. Mac\'{i}as-Virg\'{o}s, M.J. Pereira-S\'{a}ez,
A topological approach to left eigenvalues of quaternionic matrices,
Linear Multilinear Algebra 62(2) (2013)  139--158.
\bibitem{ba} A. Baker, Right eigenvalues for quaternionic matrices: a topological ap-
proach, Linear Algebra Appl. 286 (1999) 303-309.
 \bibitem{dra} T. Dray, C. A. Manogue, The octonionic eigenvalue problem, Advances
in Applied Clifford Algebras 8(2) (1998) 341-364.
\bibitem{zh} F. Zhang, Quaternions and matrices of quaternions, Linear Algebra
Appl. 251 (1997) 21-57.
\bibitem{far} D.R. Farenick, B.A.F. Pidkowich, The spectral theorem in quaternions, Linear
Algebra Appl. 371 (2003) 75-102.
\bibitem{fa} F. O. Farid,   Q.W. Wang,  F. Zhang, On the eigenvalues of quaternion matrices,
Linear Multilinear Algebra, 59(4) (2011)  451- 473.
  \bibitem{la}P. Lancaster, M. Tismenitsky,  Theory of matrices, Acad. Press., New York, 1969.
        \bibitem{liu1}X. Liu, Y. Yu, H. Wang, Determinantal representation of weighted generalized inverses, Appl. Math. Comput. 208 (2009) 556--563.
     \bibitem{liu2} X. Liu, G. Zhu, G. Zhou, Y. Yu, An Analog of the Adjugate Matrix for
the Outer Inverse ${\bf A}_{T,S}^{(2)}$, Mathematical Problems in
Engineering, Volume 2012, Article ID 591256, 14 pages,
doi:10.1155/2012/591256.
    \bibitem{ky8}I. Kyrchei, Analogs of the adjoint matrix for generalized inverses and corresponding Cramer rules, Linear  Multilinear Algebra 56(4) (2008)  453-469.
 \bibitem{ky9}I. Kyrchei,
 Explicit formulas for determinantal representations of the Drazin inverse solutions of some matrix and differential matrix equations, Appl. Math. Comput. 219 (2013) 7632-7644.
  \bibitem{ky10} I. Kyrchei, Cramer's Rule for Generalized Inverse Solutions, ,  in: Ivan I. Kyrche (Eds.), Advances in Linear Algebra Research,  Nova Sci. Publ., New York, 2015, pp. 79--132.
\bibitem{car}
Carl D. Meyer Jr.,  Limits and the index of a square matrix, SIAM J. Appl.
Math. 26(3) (1974) 506-515.


\end{thebibliography}
\end{document}